\documentclass[aap,MSNbibl,dvips]{arximspdf}
\usepackage{graphics}

\doi{10.1214/10-AAP696}
\volume{21}
\issue{1}
\pubyear{2011}
\firstpage{1}
\lastpage{32}

\makeatletter

\newtheorem{Th}{Theorem}[section]
\newproclaim{Def}[Th]{Definition}
\newtheorem{Cor}[Th]{Corollary}
\newtheorem{Lem}[Th]{Lemma}
\newtheorem{Prop}[Th]{Proposition}
\newproclaim{Rem}[Th]{Remark}

\makeatother

\begin{document}
\begin{frontmatter}

\title{Limit distributions for large P\'olya urns\protect\thanksref{T1}}
\runtitle{Limit laws for large P\'olya urns}
\thankstext{T1}{Supported by ANR-05-BLAN-0372-02 ``Structures Al\'
eatoires Discr\`etes et Algorithmes.''}

\begin{aug}
\author[A]{\fnms{Brigitte} \snm{Chauvin}\corref{}\ead[label=e1]{chauvin@math.uvsq.fr}\ead[label=u1,url]{http://www.math.uvsq.fr/\texttildelow chauvin/}},
\author[B]{\fnms{Nicolas} \snm{Pouyanne}\ead[label=e2]{pouyanne@math.uvsq.fr}\ead[label=u2,url]{http://www.math.uvsq.fr/\texttildelow pouyanne/}}
\and
\author[C]{\fnms{Reda} \snm{Sahnoun}\ead[label=e3]{sahnoun@math.uvsq.fr}}
\runauthor{B. Chauvin, N. Pouyanne and R. Sahnoun}
\affiliation{INRIA Rocquencourt and Universit\'e de Versailles,
Universit\'e de Versailles and Universit\'e de Versailles}

\address[A]{B. Chauvin\\ INRIA Rocquencourt, project Algorithms\\
Domaine de Voluceau B.P.105\\ 78153 Le Chesnay Cedex\\ France
\\
and
\\
Laboratoire de Math\'ematiques\\\quad de Versailles
\\
CNRS, UMR 8100
\\
Universit\'e de Versailles, St. Quentin
\\
45, avenue des Etats-Unis\\
78035 Versailles Cedex\\
France\\
\printead{e1}\\
\printead{u1}} 
\address[B]{N. Pouyanne\\
Laboratoire de Math\'ematiques\\\quad de Versailles\\
CNRS, UMR 8100\\
Universit\'e de Versailles, St. Quentin
\\
45, avenue des Etats-Unis\\
78035 Versailles Cedex\\
France\\
\printead{e2}\\
\printead{u2}}
\address[C]{R. Sahnoun\\
Laboratoire de Math\'ematiques de Versailles\\
CNRS, UMR 8100\\
Universit\'e de Versailles, St. Quentin\\
45, avenue des Etats-Unis\\ 78035 Versailles Cedex\\ France\\
\printead{e3}}
\end{aug}

\received{\smonth{9} \syear{2009}}
\revised{\smonth{3} \syear{2010}}

%
\begin{abstract}
We consider a two-color P\'olya urn in the case when a fixed number $S$
of balls is added at each step.
Assume it is a \emph{large} urn
that is, the second eigenvalue $m$ of the replacement matrix satisfies
$1/2<m/S\leq1$.
After $n$ drawings, the composition vector has asymptotically a first
deterministic term of order $n$ and a second random term of order $n^{m/S}$.
The object of interest is the limit distribution of this random term.

The method consists in embedding the discrete-time urn in continuous
time, getting a two-type branching process.
The dislocation equations associated with this process lead to a system
of two differential equations satisfied by the Fourier
transforms of the limit distributions.
The resolution is carried out and it turns out that the Fourier
transforms are explicitly related to Abelian integrals over the Fermat
curve of degree $m$.
The limit laws appear to constitute a new family of probability
densities supported by the whole real line.
\end{abstract}

%

\end{frontmatter}
\setcounter{footnote}{1}

\section{Introduction}\label{intro}

Consider a two-color P\'olya--Eggenberger urn random pro\-cess, with
replacement matrix
$R=\left({
a\atop c}\enskip {b\atop d }\right)$:
the urn starts with a finite number of red and black balls as initial
composition (possibly monochromatic).
At each discrete time $n$,
one draws a ball uniformly at random, notices its color,
puts it back into the urn and adds balls according to the following rule:
if the ball drawn is red, one adds $a$ red balls and $b$ black balls;
if the ball drawn is black, one adds $c$ red balls and $d$ black balls.
The integers $a,b,c,d$ are assumed to be nonnegative\footnote{One
admits classically negative values for $a$ and~$d$, together with
arithmetical conditions on $c$ and~$b$.
Nevertheless, the paper deals with so-called \emph{large} urns, for
which this never happens.}
and the urn is assumed to be \textit{balanced}, which means that the total
number of balls added at each
step is a constant $S=a+b=c+d$.
The composition vector of the urn at time $n$ is denoted by
\[
U^{\mathit{DT}}(n)=\pmatrix{
\# \mbox{ red balls at time $n$}\cr
\# \mbox{ black balls at time $n$}
}.
\]
It is a random vector and the article deals with its asymptotics when
$n$ tends to infinity.
Throughout the paper, the qualifier DT is used to refer to \emph
{discrete-time} objects while CT will refer to
\emph{continuous-time} ones.

Since P\'olya's original paper \cite{Pol}, this question has been
extensively studied so that citing all
contributions is hopeless.
The following references give however a good idea of the variety of methods:
combinatorics with many papers by  Mahmoud (see his recent book \cite
{Mah08}),
probabilistic methods by means of embedding the process in continuous
time (see  Janson \cite{Jan}),
analytic combinatorics by Flajolet {et al.} \cite{FlajDP} and a
more algebraic approach in
\cite{Pou08}.
The union of these papers is sufficiently well documented, guiding the
reader to a quasi exhaustive
bibliography.

The asymptotic behavior of $U^{\mathit{DT}}(n)$ is closely related to the
spectral decomposition of the
replacement matrix.
In case of two colors, $R$ is equivalent to
$\left({
S\atop 0}\enskip
{0\atop m}
\right)$,
where the largest eigenvalue is the balance $S$ and the smallest
eigenvalue is the integer
$m = a-c = d-b$.
We denote by $\sigma$ the ratio between the two eigenvalues:
\[
\sigma= \frac mS\leq1.
\]
It is well known that the asymptotics of the process has two different
behaviors depending on
the position of $\sigma$ with respect to the value $1/2$.
Briefly:\vspace*{1pt}

(i) when $\sigma<\frac{1}{2} $, the urn is called \emph{small}
and, except when $R$ is triangular,
the composition vector is asymptotically Gaussian\footnote{The
case $\sigma=1/2$ is similar to this one, the normalization being
$\sqrt{n\log n}$ instead of
$\sqrt n$.}:
\[
\frac{U^{\mathit{DT}}(n) -nv_1}{\sqrt n} \mathop{\hbox to 23pt{\rightarrowfill}}\limits
_{n\to\infty}\limits^{\mathcal{D}}{\mathcal N}(0,\Sigma^2),
\]
where $v_1$ is a suitable eigenvector of $^t\! R$ relative to $S$ and
$\Sigma^2$ has a simple
closed form;

(ii) when $\frac{1}{2} < \sigma< 1 $, the urn is called \emph
{large} and the composition vector
has a quite different strong asymptotic form:
%
\begin{equation}
U^{\mathit{DT}}(n) = nv_1 + n^{\sigma} W^{\mathit{DT}} v_2 +
o(n^{\sigma}),
\end{equation}
where
$v_1, v_2$ are suitable eigenvectors of $^t\! R$ relative to the
eigenvalues $S$ and $m$,
$W^{\mathit{DT}}$ is a real-valued random variable arising as the limit of a martingale,
the little $o$ being almost sure and in any $L^p, p\geq1$.
The moments of $W^{\mathit{DT}}$ can be recursively calculated but they have no
known global closed
form \cite{Pou08}.

The particular case $\sigma= 1$ is the original P\'olya urn; it
corresponds to taking $R=S\operatorname{Id}$ as replacement matrix. It
has been well
known, since Gouet \cite{Gouet}, that the composition vector admits an
almost sure asymptotics of order one: $U^{\mathit{DT}}(n) = n D +o(n)$ where the
random vector $D$ has a Dirichlet density (explicitly given in~\cite{Gouet}).

In the present article, the object of interest is the distribution of
$W^{\mathit{DT}}$ for large urns.

The first step consists in classically embedding the discrete-time
process $(U^{\mathit{DT}}(n))_{n\geq0}$
into a continuous-time Markov process $(U^{\mathit{CT}}(t))_{t\geq0}$, by
equipping each ball with an exponential
clock.
At any $n$th jump time $\tau_n$, the continuous-time process
$U^{\mathit{CT}}(\tau_n)$ has the same distribution as $U^{\mathit{DT}}(n)$.
This connection between both processes is the key point, allowing us to
work on the continuous-time process, where independence properties have
been gained.

In Theorem \ref{continuousurn}, we show that, in the case of large
urns, the continuous-time process
satisfies, when $t$ tends to infinity, the following asymptotics,
%
\begin{equation}
U^{\mathit{CT}}(t) = e^{S t}\xi v_1 \bigl(1+o(1)\bigr)+ e^{m t} W^{\mathit{CT}} v_2\bigl(1+o(1)\bigr),
\end{equation}
where $v_1$ and $v_2$ are the same vectors as above, $\xi$ and
$W^{\mathit{CT}}$ are real-valued random variables that arise as limits of
martingales, with $o(\cdot)$ meaning ``almost sure and in any $L^p, p\geq1$.''
Moreover, we prove that $\xi$ is Gamma-distributed.
These results are based on the spectral decomposition of the
infinitesimal generator of the continuous-time process on spaces of
two-variable polynomials.

Thanks to the embedding connection, the two random variables $W^{\mathit{DT}}$
and $W^{\mathit{CT}}$ are connected (Theorem~\ref{thmartingaleconnection}):
\[
W^{\mathit{CT}} = \xi^{\sigma}  W^{\mathit{DT}}   \qquad\mbox{a.s.},
\]
$\xi$ and $W^{\mathit{DT}}$ being independent.
Since $\xi^{\sigma}$ is invertible,\footnote{A
probability distribution $A$ is called \emph{invertible} when, for
any probability distributions $A$ and
$B$, the equation $AX=B$ admits a unique solution $X$ independent of
$A$, see, for instance, Chaumont and Yor
 \cite{ChaumontYor}.
The invertibility of any power of a Gamma distribution can be shown by
elementary considerations
on Fourier transforms.} the attention is focused on the determination of the distribution of
$W^{\mathit{CT}}$.

Because of the nonnegativity of $R$ entries, $(U^{\mathit{CT}}(t))_{t\geq0}$ is
a two-type
branching process, visualized as a random tree:
the branching property gives rise to dislocation equations on $U^{\mathit{CT}}(t)$.
If one denotes by $\mathcal{F}$ ({resp.}, $\mathcal{G}$) the
characteristic function of $W^{\mathit{CT}}$ starting
from one red ball and no black ball ({resp.}, no red ball, one
black ball), the independence of the
subtrees in the branching process implies that the characteristic
function of \emph{any} $W^{\mathit{CT}}$
starting from $\alpha$ red balls and $\beta$ black balls is the
product $\mathcal{F}^{\alpha}\mathcal{G}^{\beta}$.
Furthermore, the dislocation equations on $U^{\mathit{CT}}(t)$ lead to the
following differential system
%
\begin{equation}
\label{premiersystemeFourier}
\cases{
\mathcal{F}(x)+mx\mathcal{F}'(x)=\mathcal{F}(x)^{a+1}\mathcal
{G}(x)^b,
\cr
\mathcal{G}(x)+mx\mathcal{G}'(x)=\mathcal{F}(x)^{c}\mathcal{G}(x)^{d+1},
}
\end{equation}
together with suitable boundary conditions.
Notice that the corresponding exponential moment generating series
(Laplace series) are also
solutions of (\ref{premiersystemeFourier}), but their radius of
convergence is equal to 0.
This is detailed in Section \ref{serieLaplace}.

The solution of system (\ref{premiersystemeFourier}) is obtained in
Section \ref{resolution}, where it is shown that $\mathcal{F}$ and
$\mathcal{G}$ can be made explicit in terms of inverse functions of Abelian
integrals over the Fermat curve of degree $m$. Indeed, for any complex
$z$ in a suitable open subset of $\mathbb{C}$, let
\[
I_{m,S,b}(z)=\int_{[z,z\infty)}(1+u^m)^{b/m}\frac
{du}{u^{S+1}},
\]
where $[z,z\infty)$ denotes the ray $\{ tz, t\geq1\}$.
The function $I_{m,S,b}$ defines a conformal mapping   on the
open sector
$
\mathcal{V}_m=\{ z\neq0, 0<\arg(z)<\pi/m\}
$.
If $J_{m,S,b}$ denotes the holomorphic function, defined on the lower
half-plane as left inverse function of $I_{m,S,b}$ and extended to the
slit plane by conjugacy, the closed
form of $\mathcal{F}$ and $\mathcal{G}$ are given in the following result.

\begin{Th*}
For any $x>0$,\vspace*{-2pt}
\[
\cases{
\mathcal{F}(x)=Kx^{-1/m}J_{m,S,b}\biggl( C_0+\dfrac
{K^S}Sx^{-
S/m}\biggr),
\cr
\mathcal{G} (x)=Kx^{-1/m}J_{m,S,c}\biggl( C_0+\dfrac
{K^S}Sx^{- S/m}\biggr) ,
}\vspace*{-2pt}
\]
where $K\in\mathbb{C}$ and $C_0<0$ are explicit constants.
\end{Th*}

For precise statements and proofs, see Section \ref{calculdeF} and Theorem
\ref{characteristicFunctions}.

The solution of system (\ref{premiersystemeFourier}) is effected by a
ramified change of variable and
functions, leading to the following monomial system:
%
\begin{equation}
\label{monom}
\cases{
f'=f^{a+1}g^b,
\cr
g'=f^cg^{d+1}.
}
\end{equation}
This remarkable fact is evocative of the case of \emph{small} urns and
discrete time, as considered in a beautiful study of Flajolet {et
al.} \cite{Flaj}. The method of \cite{Flaj} leads directly to the
same system (\ref{monom}) on generating functions.
The assumption $\sigma<1/2$, when expressed in term of the four
parameters $a$, $b$, $c$ and $d$, does not fundamentally affect the
system but requires completely different analytic
considerations.

The limit laws of the $W^{\mathit{CT}}$ appear to constitute a new family of
probability distributions, indexed by
three parameters $S,m,b$ subject to assumptions (\ref{hypSmb}) and by
initial conditions
$\alpha, \beta$.
We prove in Section \ref{densiteWCT} that they admit densities that
can be expressed by means of the inverse Fourier transforms of their
characteristic function derivatives.
Furthermore, the laws are infinitely divisible and their support is the
whole real line, the radius of
convergence of their exponential moment generating series being equal
to $0$.

Many questions remain open.
For instance, are these distributions characterized by their moments?
What is the precise asymptotics of their densities at infinity (tails)?
It is shown in \cite{Sahnoun2} that, for triangular and nondiagonal
replacement matrices, the
discrete-time limit law $W^{\mathit{DT}}$ is never infinitely divisible;
does this situation persist in the present nontriangular case?



\section{The model}\label{model}

\subsection{Definition of the process}\label{sec2.1}

Let $a$, $b$, $c$ and $d$ be nonnegative integers such that $a+b=c+d
=:S$ and $R$ be the matrix
%
\begin{equation}
\label{matriceUrne}
R:=
\pmatrix{
a&b\cr
c&d
}
=
\pmatrix{
a&S-a\cr
S-d&d
}.
\end{equation}
The discrete-time P\'olya--Eggenberger urn process associated with the
replacement matrix $R$,
which has been informally described in the \hyperref[intro]{Introduction}, is the Markov
chain $(U^{\mathit{DT}}(n), n\in\mathbb{N})$, having $\mathbb{N}^2\setminus
\{ (0,0)\}$ as
state space and
%
\begin{equation}
\label{probasTransition}
\frac x{x+y}\delta_{(x+a,y+b)}+\frac y{x+y}\delta_{(x+c,y+d)}
\end{equation}
as transition probability at any nonzero point $(x,y)\in\mathbb
{N}^2$. In
this formula, $\delta_v$ denotes Dirac point mass at $v\in\mathbb{N}^2$.
This means that $(U^{\mathit{DT}}(n), n\in\mathbb{N})$ is a random walk in
$\mathbb{N}^2\setminus\{ (0,0)\}$
(or in the two-dimensional one-column nonzero matrices with nonnegative
integer entries, we will use
both notations)
recursively defined
by the conditional probabilities
\[
\cases{
\mathrm{P}
\biggl( U^{\mathit{DT}}(n+1) = U^{\mathit{DT}}(n) + \pmatrix{ a\cr b }
| U^{\mathit{DT}}(n)=\pmatrix{ x\cr y }\biggr)
=\dfrac{x}{x+y},
 \cr
\mathrm{P}
\biggl( U^{\mathit{DT}}(n+1) = U^{\mathit{DT}}(n) + \pmatrix{ c\cr d }
| U^{\mathit{DT}}(n)=\pmatrix{ x\cr y }\biggr)
=\dfrac{y}{x+y}.
}
\]
In the sequel,
\[
\bigl( U_{(\alpha,\beta)}^{\mathit{DT}}(n),n\geq0\bigr)
\]
will denote the process starting from the nonzero vector $(\alpha
,\beta)$ and
\[
u:=\alpha+\beta
\]
will denote the total number of balls at time $0$.
Notice that the balance property $S=a+b=c+d$ implies that the total
number of balls at time $n$,
when $U^{\mathit{DT}}(n)=(x,y)$, is the (nonrandom) number $x+y=u+nS$.

Denote by $w_1={ a\choose b }$ and $w_2={ c\choose d }$ the
increment vectors of the walk.
The transition operator $\Phi$ is defined, for any function $f$ on
$\mathbb{N}^2$ and for any $v={ x\choose y }$, by
\[
\Phi(f)(v) = x[ f(v+w_1) - f(v)] + y[ f(v+w_2) -
f(v)] .
\]
Conditioning on $({\mathcal F}_n, n\geq0)$, which is the filtration
associated with the process $(U^{\mathit{DT}}(n) , n\geq0)$, one gets
\[
\mathbb{E}^{{\mathcal F}_n}\bigl( f\bigl( U^{\mathit{DT}}(n+1)\bigr)\bigr) = \biggl( I + \frac
{\Phi
}{u+nS}\biggr) (f) ( U^{\mathit{DT}}(n) ).
\]
In particular,
%
\begin{equation}
\label{mart}
\mathbb{E}^{{\mathcal F}_n}\bigl( U^{\mathit{DT}}(n+1)\bigr) = \biggl( I+\frac
{A}{u+nS}\biggr) U^{\mathit{DT}}(n),
\end{equation}
where $I$ denotes the two-dimensional identity matrix and
\[
A := {}^t\hspace*{-2pt}R
=\pmatrix{
a&c\cr b&d
}.
\]

\subsection{Asymptotics of the discrete-time process $U^{\mathit{DT}}(n)$}

As mentioned in Section \ref{intro} and briefly recalled hereunder, a
discrete-time P\'olya--Eggenberger
urn process has two different kinds of asymptotics depending on the
ratio of the eigenvalues of its replacement matrix $R$.
With our notation, these eigenvalues are $S$ and\vspace*{-1pt}
\[
m:=a-c=d-b.
\]
Let us denote by $u_1$ and $u_2$ the two following linear
eigenforms\footnote{An eigenform of an endomorphism $f$ is an
eigenvector of $^t\hspace*{-2pt} f$; some authors call these linear forms left
eigenvectors of $f$, referring to matrix operations.} of $A$,
respectively associated with the eigenvalues $S$ and $m$, which means
that $u_1\circ A=Su_1$ and $u_2\circ A=mu_2$:\vspace*{-1pt}
%
\begin{equation}\label{u}
u_1(x,y) = \frac{1}{S}(x+y), \qquad
u_2(x,y)= \frac{1}{S}(bx-cy),
\end{equation}
and denote by $(v_1,v_2)$ the dual basis of $(u_1,u_2)$:\vspace*{-1pt}
%
\begin{equation}\label{v}
v_1 = \frac S{(b+c)} \pmatrix{
c\cr
b
},
\qquad v_2 = \frac S{(b+c)} \pmatrix{
1\cr
-1
}.
\end{equation}
The vectors $v_k$ are eigenvectors of $A$ and the projections on the
eigenlines are $u_1v_1$
and $u_2v_2$.

For any positive real $x$ and any nonnegative integer $n$, if one
denotes by $\gamma_{x,n}$ the
polynomial\vspace*{-1pt}
\[
\gamma_{x,n} (t):= \prod_{k=0}^{n-1} \biggl(1 + \frac t{x+k}\biggr),
\]
the matrix $\gamma_{\frac mS,n}(\frac AS)$ in nonsingular and it is
immediate from (\ref{mart}) that\vspace*{-1pt}
\[
\gamma_{ \frac mS,n}\biggl(\frac AS\biggr)^{-1}  U^{\mathit{DT}}(n)
\]
is a (vector-valued) martingale.

As indicated in the \hyperref[intro]{Introduction}, the ratio of $R$ eigenvalues is
denoted by\vspace*{-1pt}
\[
\sigma:=\frac mS\leq1.
\]
The case of \emph{small} urn processes ({i.e.,} when $\sigma\leq1/2$) has been well studied;
in this case, when $R$ is not triangular, the random vector admits a
Gaussian central limit theorem
(see Janson  \cite{Jan}).
Triangular replacement matrices impose a particular treatment and lead
most often to a nonnormal
second-order limit (see Janson  \cite{JanTrig} or \cite{Sahnoun2}).

Our subject of interest is the case of \emph{large} urns, that is,
when $\sigma>1/2$.
In this case, $\frac1{S} U^{\mathit{DT}}(n)$ is a large P\'olya process with
replacement matrix $\frac1{S} R$
in the sense of \cite{Pou08}.
As a matter of consequence, the projections of the above vector-valued
martingale on the eigenlines
of $A$, which are of course also martingales, converge in any ${\rm
L}^p$, $p\geq1$ (and a.s.).
In particular (second projection),
\[
M^{\mathit{DT}}(n):= \frac{1}{\gamma_{\frac uS,n}(\sigma) } u_2(U^{\mathit{DT}}(n))
\]
is a convergent martingale;
since $\gamma_{u,n}(\sigma) = n^{\sigma} \frac{\Gamma(u)}{\Gamma
(u+\sigma)}(1+o(1))$, denoting by
%
\begin{equation}\label{martingalediscrete}
W^{\mathit{DT}} := \lim_{n\rightarrow+\infty} \frac{1}{n^{\sigma}} u_2(U^{\mathit{DT}}(n)),
\end{equation}
a slight adaptation of \cite{Pou08} leads to the following theorem.
Note that this theorem was essentially proven by Athreya and Karlin
\cite{AK} and Janson \cite{Jan}
for random replacement matrices.
The convergence in $L^p$-spaces when $R$ is nonrandom is shown by the
indicated adaptation
of \cite{Pou08}.

\begin{Th}
\label{asymptotiquediscrete}
Suppose that $\sigma\in\,]1/2,1[$.
Then, as $n$ tends to infinity,
\[
U^{\mathit{DT}}(n) = nv_1 + n^{\sigma} W^{\mathit{DT}} v_2 + o( n^{\sigma} ),
\]
where $v_1$ and $v_2$, defined in \textup{(\ref{v})}, are eigenvectors
associated with the eigenvalues $S$ and $m$; $W^{\mathit{DT}}$ is defined by
(\ref{martingalediscrete}); $o( \cdot)$ means a.s. and in any $L^p, p\geq1$.
\end{Th}

\subsection{Parametrization and hypotheses}
\label{parametrage}

The subject of the paper is $W^{\mathit{DT}}$ distribution in Theorem \ref
{asymptotiquediscrete} so that
the P\'olya urn process will be supposed large.
Furthermore, the replacement matrix $R$ is supposed to be not
triangular because this case has
to be treated separately with regard to its limit law, as attested by
Janson \cite{JanTrig}, Flajolet et al. \cite{FlajDP}, \cite{Sahnoun2} and the present paper.

\begin{figure}[t]

\includegraphics{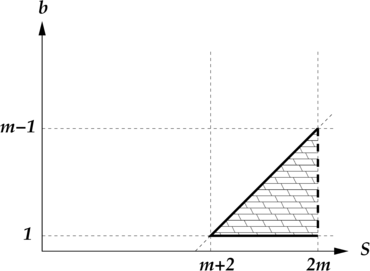}

\caption{Parameters $(b,S)$ for a given $m$.}\label{fig1}
\end{figure}

Under these conditions, the assumptions on the replacement matrix
$R=\left(
{a\atop c}\enskip {b\atop d
}
\right)$
are: $a+b=c+d=S$ (balance condition), $S/2<m=a-c=d-b<S$ (large urn) and
$b,c\geq1$
(not triangular).
Because of the balance condition, the parametrization of P\'olya urns
is governed by three free
parameters.
The computation of Fourier transforms will show in Section \ref
{calculdeF} that a natural
choice of free parameters is $(m,S,b)$.
The assumption ``large and nontriangular'' is equivalent, in terms of
these data, to the following:
\[
R=\pmatrix{a&b\cr c&d
}
=\pmatrix{S-b&b\cr S-m-b&m+b
}
\]
with
%
\begin{equation}
\label{hypSmb}
\cases{
m+2\leq S\leq2m-1,
\cr
1\leq b\leq S-m-1.
}
\end{equation}
Note that these inequalities imply $S\geq5$ and $m\geq3$ and that,
for a given $m$, the point $(m,b)$
belongs to a triangle as represented in Figure~\ref{fig1}.

For small values of $S$, large urn processes have the following
possible replacement matrices:
for $S\in\{1,2\}$, only $R=S\operatorname{Id}_2$ defines a large urn;
for $S\in\{3,4\}$, all large urns have triangular matrices.
For $S\in\{5,6\}$, only
$R=
\left(
{
S-1\atop 1} \enskip {1\atop S-1
}
\right)
$
defines a nontriangular large urn.
For $S=7$, apart from triangular or symmetric matrices, there are only
two replacement matrices that define large urns: $
\left(
{
6\atop 2}\enskip{1\atop 5
}
\right)
$
and the other one derived from it by permutation of coordinates.


\section{Embedding in continuous time and martingale connection}\label{method}

\subsection{Embedding}\label{embeddingSection}

The idea of embedding discrete urn models in continuous-time branching
processes goes back
at least to Athreya and Karlin \cite{AK}.
A description is given in Athreya and Ney  \cite{AN}, Section 9.
The method has been recently revisited and developed by Janson \cite{Jan}.

We define the continuous-time Markov branching process $(U^{\mathit{CT}}(t),
t\geq0)$ as
being the embedded process of $(U^{\mathit{DT}}(n), n\geq0)$.
Following, for instance, Bertoin \cite{Bertoin}, Section 1.1, this
means that it is defined as the continuous-time Markov chain having as
jump rate, at any nonzero point $(x,y)\in\mathbb{N}^2$, the finite
measure given by the transition probability of the discrete-time process
[formula (\ref{probasTransition})].
One gets this way a branching process having the following dynamical
description in terms of red and black balls.
In the urn, at any moment, each ball is equipped with an $\mathit{\mathcal{E}xp}(1)$-distributed\footnote{
For any positive real $a$, $\operatorname{\mathit{\mathcal{E}xp}}(a)$ denotes the exponential
distribution with parameter $a$.} random clock, all the clocks being
independent.
When the clock of a red ball rings, $a$ red balls and $b$ black balls
are added in the urn;
when the ringing clock belongs to a black ball, one adds $c$ red balls
and $d$ black balls, so that
the replacement rules are the same as in the discrete-time urn process.

The successive jumping times of $(U^{\mathit{CT}}(t), t\geq0)$, will be denoted by
\[
0=\tau_0 < \tau_1 < \cdots< \tau_n < \cdots.
\]
The $n$th jumping time is the time of the $n$th dislocation of the
branching process.
The process is thus constant on any interval $[\tau_n,\tau_{n+1}[$.

In the sequel,
\[
\bigl( U_{(\alpha,\beta)}^{\mathit{CT}}(t),t\geq0\bigr)
\]
will denote the process starting from the nonzero vector $(\alpha
,\beta)$.
Thus, for any initial condition $(\alpha, \beta)$, for any $t\geq0$,
\[
U_{(\alpha,\beta)}^{\mathit{CT}}(t)= U_{(\alpha,\beta)}^{\mathit{DT}}(
a(t)),
\]
where
\[
a(t):= \inf\{ n\geq0, \tau_n \geq t \}.
\]


\begin{Lem}
\label{embed}
(i) for $n\geq0$, the distribution of $\tau_{n+1} - \tau_{n}$ is
$\mathit{\mathcal{E} xp} (u+Sn)$ where $u$
denotes the total number of balls at time $0$;\vspace*{-6pt}
\begin{longlist}[(iii)]
\item[(ii)] the processes $(\tau_n)_{n\geq0}$ and $(U^{\mathit{CT}}(\tau_n))_{n\geq
0}$ are independent;

\item[(iii)] the processes $(U^{\mathit{CT}}(\tau_n))_{n\geq0}$ and $(U^{\mathit{DT}}(n))_{n\geq
0}$ have the same distribution.
\end{longlist}
\end{Lem}

\begin{pf}
The total number of balls at time $t\in[\tau_n,\tau_{n+1}[$ is $u+Sn$.
Therefore, (i) is a consequence of the fact that the minimum of $k$
independent random variables
$\mathit{\mathcal{E} xp} (1)$-distributed is $\mathit{\mathcal{E} xp} (k)$-distributed.
(ii) is the classical independence between the jump chain and the jump
times in such Markov processes.
The initial states and evolution rules of both Markov chains in
discrete time and in continuous time are
the same ones, so that (iii) holds.
\end{pf}

\textit{Convention}: From now on, thanks to (iii) of Lemma
\ref{embed}, we will classically
consider that the discrete-time process and the continuous-time process
are built on the \emph{same}
probability space on which
%
\begin{equation}
\label{embedding}
(U^{\mathit{CT}}(\tau_n))_{n\geq0} = (U^{\mathit{DT}}(n))_{n\geq0}  \qquad \mbox{a.s.}
\end{equation}

\subsection{Asymptotics of the continuous-time process $U^{\mathit{CT}}(t)$}
\label{asymptCT}

Let $v_1$ and $v_2$ the linearly independent eigenvectors of $A$
defined by (\ref{v}).
In the case of large urns, the asymptotics of the continuous-time
process $(U^{\mathit{CT}}(t))_{t\geq0}$
is given in the following theorem.

\begin{Th}[(Asymptotics of continuous-time process)]\label{continuousurn}
When $t$ tends to infinity,
%
\begin{equation}
\label{asymptotiqueCT}
U^{\mathit{CT}}(t)=e^{St} \xi v_1\bigl( 1+o(1)\bigr) +e^{mt}W^{\mathit{CT}}v_2\bigl(
1+o(1)\bigr) ,
\end{equation}
where $\xi$ and $W^{\mathit{CT}}$ are real-valued random variables, the little
$o$ being almost sure and
in any ${\rm L}^p$-space, $p\geq1$.
Furthermore, $\xi$ is  $\operatorname{Gamma}(u /S)$-distributed, where $u
=\alpha+\beta$ is the total
number of balls at time $0$.
\end{Th}

\begin{Rem}
Another formulation of this Theorem is:
in the basis $(v_1,v_2)$, the coordinate of $U^{\mathit{CT}}(t)$ along $v_1$ has
$e^{St} \xi$ as its dominant term while the coordinate of $U^{\mathit{CT}}(t)$
along $v_2$ has $e^{mt} W^{\mathit{CT}}$ as its dominant term.
\end{Rem}

\begin{pf*}{Proof of Theorem~\ref{continuousurn}}
The embedding in continuous time has been studied in
Athreya and Karlin \cite{AK} and in
Janson \cite{Jan}.
It has become classical  that the process
\[
( e^{-tA}U^{\mathit{CT}}(t)) _{t\geq0}
\]
is a vector-valued martingale and that, in the case of large urns
($\sigma>1/2$), this martingale is
bounded in ${\rm L}^2$, thus converges.
Its projections on the eigenlines $\mathbb{R}v_1$ and $\mathbb
{R}v_2$, \emph
{\textup{that is},} respectively,
\[
e^{-St}u_1( U^{\mathit{CT}}(t))
\quad\mbox{and}\quad
e^{-mt}u_2( U^{\mathit{CT}}(t))
\]
are also ${\rm L}^2$-convergent real valued martingales, thus converge
almost surely.
Their respective limits are named $\xi$ and $W^{\mathit{CT}}$.
What still has to be proved is that these martingales converge in fact
in any ${\rm L}^p$, $p\geq1$.
The identification of $\xi$ distribution will be a consequence of this proof.

The infinitesimal generator of the Markov process $(U^{\mathit{CT}}(t))_t$ is
the finite-difference operator
\[
\Phi(f)(x,y)=x\{ f(x+a,y+b)-f(x,y)\}
+y\{ f(x+c,y+d)-f(x,y)\},
\]
defined for any (measurable) function $f$ and any $(x,y)\in\mathbb{R}^2$.
For a very synthetic reference on semi-groups of Markov continuous-time
processes,
one can refer to Bertoin \cite{Bertoin}, Chapter 1.
This operator $\Phi$ acts on two-variable  polynomials.
This action has been studied in detail in \cite{Pou08} in a more
general framework.
More precisely, for any integer $d\geq1$, the operator $\Phi$ acts on
the finite-dimensional space
of polynomials of degree less than $d$, so that, for any two-variable
polynomial $P$ and for any
$t\geq0$,
%
\begin{equation}
\label{generateurInfinitesimal}
\mathbb{E}( P(U^{\mathit{CT}}(t))) =\exp( t\Phi
)(
P) ( (U^{\mathit{CT}}(0))),
\end{equation}
where, in this formula, $\Phi$ denotes the restriction of $\Phi$
itself on any finite-dimensional
polynomials space containing $P$.
The properties of $\Phi$ listed in the following lemma are proven in
\cite{Pou08}
and will be used here.

\begin{Lem}
\label{lemmePhi}
There exists a unique family of polynomials $Q_{p,q}\in\mathbb{R}[x,y]$,
$p,q$ nonnegative integers, called \emph{reduced polynomials},
such that:
\begin{longlist}[(1)]
\item[(1)]
$Q_{0,0}=1$, $Q_{1,0}=u_1$ and $Q_{0,1}=u_2$
[see (\ref{u}) for a definition of eigenforms $u_1$ and $u_2$];

\item[(2)]
$\Phi(Q_{p,q})=(pS+qm)Q_{p,q}$ for all nonnegative integers $p,q$;

\item[(3)]
$u_1^pu_2^q-Q_{p,q}\in\operatorname{Span}\{ Q_{p',q'}, p'S+q'm<pS+qm\}
$ for all
nonnegative integers $p,q$.
\end{longlist}
\end{Lem}

Note that the reduced polynomial $Q_{p,q}$ is in fact the projection of
$u_1^pu_2^q$ on a suitable
characteristic subspace of $\Phi$ restriction to some finite-dimensional polynomial space, and that
this spectral decomposition of $\Phi$ on polynomials has a
particularly simple form (it is
diagonalizable) because the urn is large and two-dimensional.
See \cite{Pou08} for more details.

Formula (\ref{generateurInfinitesimal}) and property {(ii)} of
Lemma \ref{lemmePhi}
lead to
\[
\forall(p,q)\in\mathbb{Z}_{\geq0}^2\qquad
\mathbb{E}( Q_{p,q}( U^{\mathit{CT}}(t)))
=e^{t(pS+qm)}
Q_{p,q}(U^{\mathit{CT}}(0)).
\]
This implies straightforwardly, with {(iii)} of Lemma \ref
{lemmePhi}, that, for any $(p,q)$,
%
\begin{equation}
\label{momentsReduits}
\mathbb{E}( u_1^pu_2^q( U^{\mathit{CT}}(t)))
=e^{t(pS+qm)} Q_{p,q}(U^{\mathit{CT}}(0))
+o\bigl( e^{t(pS+qm)}\bigr).
\end{equation}
In particular, the martingales $e^{-St}u_1( U^{\mathit{CT}}(t))$
and $e^{-mt}u_2( U^{\mathit{CT}}(t))$
are ${\rm L}^p$-bounded for any $p\geq1$ and their respective limits,
namely $\xi$ and $W^{\mathit{CT}}$
satisfy, for any nonnegative integer $p$,
%
\begin{equation}
\label{momentsXiW}
\mathbb{E}\xi^p=Q_{p,0}( U^{\mathit{CT}}(0))
 \quad\mbox{and}\quad
\mathbb{E}( W^{\mathit{CT}}) ^p=Q_{0,p}( U^{\mathit{CT}}(0)).
\end{equation}
The convergence part of the theorem follows now from the spectral
decomposition of $A$:
for any $t\geq0$,
\[
U^{\mathit{CT}}(t)=u_1( U^{\mathit{CT}}(t))\cdot v_1+u_2( U^{\mathit{CT}}(t))\cdot v_2.
\]
Besides, it is proven in \cite{Pou08}, or one can check it after an
easy computation, that the reduced
polynomials corresponding to the powers of $u_1$ have the following
closed form expression
\[
Q_{p,0}=u_1( u_1+1)( u_1+2)\cdots(
u_1+p-1).
\]
Thanks to formula (\ref{momentsXiW}), this shows that the $p$th moment
of $\xi$ is, for any integer
$p\geq0$,
\[
\mathbb{E}\xi^p=\frac u S\biggl( \frac u S +1\biggr)\biggl( \frac u S
+2\biggr)\cdots\biggl( \frac u S +p-1\biggr)
=\frac{\Gamma( \frac{u}{ S}+p)}{\Gamma(  \frac{u}{
S})},
\]
where $u$ is the total number of balls at time $0$ [remember that
$u_1(U^{\mathit{CT}}(0))=u /S$,
see (\ref{u})].
One identifies this way the required Gamma$(u /S)$ distribution,
characterized by its moments.
\end{pf*}

\begin{Rem}
Notice that the distribution of $\xi$ has been given by Janson \cite
{Jan} calculating first the distribution of $u_1(U^{\mathit{CT}}(t))$for every $t$:
\[
u_1(U^{\mathit{CT}}(t)) = \frac u S + Z(t),
\]
where $Z(t)$ is a negative binomial distribution.
\end{Rem}

\begin{Rem}
Reduced polynomials $Q_{0,p}$ do not have a known closed form, so that
reproducing the above method in order to compute the moments of
$W^{\mathit{CT}}$ fails.
\end{Rem}

\begin{Rem}
It follows from the proof that the real-valued random variables $\xi$
and $W^{\mathit{CT}}$ are respective limits
of the martingales
\begin{eqnarray*}
\xi&=&\lim_{t\to+\infty}e^{-St}u_1( U^{\mathit{CT}}(t)),
\\
W^{\mathit{CT}} &=&\lim_{t\to+\infty}e^{-mt}u_2( U^{\mathit{CT}}(t)).
\end{eqnarray*}
They are not independent and their joint moments are computed from
formula~(\ref{momentsReduits}):
for any nonnegative integers $p,q$,
\[
E[ ( \xi) ^p( W^{\mathit{CT}}) ^q
]=Q_{p,q}( U^{\mathit{CT}}(0)).
\]
For example, their respective means are $E\xi=u_1(U^{\mathit{CT}}(0))=\frac
1S(\alpha+\beta)$ and
$EW^{\mathit{CT}}=u_2(U^{\mathit{CT}}(0))=\frac1S(b\alpha-c\beta)$, whereas
\[
E( \xi W^{\mathit{CT}}) =\frac{(\alpha+\beta+m)(b\alpha-c\beta)}{S^2}
\]
as can be shown by computation of the
reduced polynomial $Q_{1,1}=(u_1+\sigma)u_2$ (one can directly check
that this polynomial is
an eigenvector of $\Phi$, associated with the eigenvalue $S+m$).
\end{Rem}

\begin{Rem}
When the urn is small ($\sigma<1/2$), the same method shows that the
result on the first projection is
still valid:
the martingale\break $(e^{-St}u_1(U^{\mathit{CT}}(t)))_t$ converges in any ${\rm L}^p$
($p\geq1$) to a
$\operatorname{Gamma}(u /S)$ distributed random variable.
On the contrary, the martingale $(e^{-mt}u_2(U^{\mathit{CT}}(t)))_t$ diverges
and it is shown in
Janson \cite{Jan} that
the second projection satisfies a central limit theorem:
when $\sigma=\frac mS<1/2$,
\[
e^{-  \frac S2 t}u_2( U^{\mathit{CT}}(t)) \mathop{\hbox to 23pt{\rightarrowfill}
}\limits_{t\to+\infty}\limits^{\mathcal{D}}\mathcal{N},
\]
where $\mathcal{N}$ is a normal distribution.
In the case $\sigma=1/2$, the normalization must be modified and one
gets the convergence in law of
$\sqrt te^{-St/2}u_2(U^{\mathit{CT}}(t))$ to a normal distribution.
\end{Rem}

\begin{Rem}
The distributions of the $W^{\mathit{CT}}$ are infinitely divisible, because
they are limits of infinitely divisible ones, obtained by scaling and
projection of infinitely divisible ones. Indeed, in finite time, the
distributions of the $U_{(\alpha, \beta)}^{\mathit{CT}}(t) $ are infinitely
divisible. It has been said by Janson \cite{Jan}, proof of Theorem
3.9.
With our notations, especially the one of (\ref{dislocationVectorielle}), it relies on the fact that
\[
U_{(\alpha, \beta)}^{\mathit{CT}}(t) \stackrel{\mathcal{L}}{=}[n] U_{(
{\alpha}/n,
{\beta}/n)}^{\mathit{CT}}(t),
\]
where a continuous-time branching process (with the same branching
dynamics as before), starting from real (nonnecessary integer)
conditions, is suitably defined.
\end{Rem}

\subsection{DT and CT connections}
\label{connections}

Apply the first projection to the embedding principle (\ref{embedding}):
\[
u_1 ( U^{\mathit{CT}}(\tau_n)) = u_1 (U^{\mathit{DT}}(n) )  \qquad \mbox{a.s}.
\]
By definition (\ref{u}) of $u_1$, this number is $\frac{1}S$ times
the number of balls in the urn at time
$n$, which equals $\frac{1}S (u +Sn) = n (1+o(1))$.
Since stopping times $\tau_n$ tend to $+\infty$, renormalizing by
$e^{-S\tau_n}$ and applying the
convergence result of Section \ref{asymptCT} leads to
%
\begin{equation}\label{limxi}
\xi= \lim_{ n\rightarrow+\infty} n e^{ -S\tau_n}.
\end{equation}
Apply now the second projection to the embedding principle (\ref{embedding}):
\[
u_2 ( U^{\mathit{CT}}(\tau_n)) = u_2 (U^{\mathit{DT}}(n) ) \qquad \mbox{a.s}.
\]
Renormalizing by $e^{-m\tau_n}$ implies that
%
\[
e^{ -m \tau_n} u_2 ( U^{\mathit{CT}}(\tau_n)) = W^{\mathit{CT}}(\tau_n ) = e^{ -m \tau
_n} \gamma_{\frac{u}{S},n}(\sigma) M^{\mathit{DT}}(n)\qquad \mbox{a.s}.
\]
which is a ``martingale connection'' in finite time.

Thanks to (\ref{limxi}) and Theorem \ref{asymptCT},
passing to the limit $n\to\infty$ leads to the following theorem,
already mentioned in Janson \cite{Jan}
in a more general framework.
Note that the independence between $\xi$ and $W^{\mathit{DT}}$ comes from Lemma
\ref{embed}(ii).

\begin{Th}[(Martingale connection)]
\label{thmartingaleconnection}
%
\begin{equation}
\label{martingaleconnection}
W^{\mathit{CT}} = \xi^{\sigma}  W^{\mathit{DT}} \qquad \mbox{a.s}.
\end{equation}
$\xi$ and $W^{\mathit{DT}}$ being independent.
\end{Th}

The distribution of $\xi^{\sigma}$ is invertible (see footnote in the
\hyperref[intro]{Introduction}),
so that any information on $W^{\mathit{CT}}$ can be pulled back to $W^{\mathit{DT}}$
thanks to connection
(\ref{martingaleconnection}).


\section{Dislocation equations for continuous urns}\label{dislocation}

\subsection{Vectorial finite time dislocation equations}
\label{dislocationvect}

By embedding in continuous time, the previous section provided a
branching process $(U_{(\alpha, \beta)}^{\mathit{CT}}(t), t\geq0)$. The
independence properties of this process imply that it is equal to the
sum of $\alpha$ copies of $U_{(1,0)}^{\mathit{CT}}(t)$ (the process starting
from one red ball) and $\beta$ copies of $U_{(0,1)}^{\mathit{CT}}(t)$ (the
process starting from one black ball). We are led to study these two
$\mathbb{R}^2$-valued processes.

Let us now apply the strong Markov branching property to these
processes: let us denote by $\tau$ the first splitting time for any of
these processes (they have the same $\mathit{\mathcal{E} xp} (1)$
distribution). We
get the following vectorial finite time dislocation equations:
%
\begin{equation}
\label{dislocationVectorielle}
\forall t>\tau
\cases{
U_{(1,0)}^{\mathit{CT}}(t)\stackrel{\mathcal{L}}{=}[a+1]U_{(1,0)}^{\mathit{CT}}(t-\tau
)+[b]U_{(0,1)}^{\mathit{CT}}(t-\tau),
 \cr
U_{(0,1)}^{\mathit{CT}}(t)\stackrel{\mathcal{L}}{=}[c]U_{(1,0)}^{\mathit{CT}}(t-\tau
)+[d+1]U_{(0,1)}^{\mathit{CT}}(t-\tau),
}
\end{equation}
where the notation $[n]X+[m]Y$ stands for the sum of $n$ copies of the
random variable $X$ and $m$ copies of the random variable $Y$ ($n$ and
$m$ are nonnegative integers).

\begin{Rem}
The above equations could be written with a.s. equalities. Taking a
probability space of trees is more convenient. The price to pay is just
to write the different processes for each subtree with different
indexes and to distinguish the two splitting times for the two starting
situations.
\end{Rem}

\subsection{Limit dislocation equations}\label{dislocationlimite}

Remember that $( e^{-mt}u_2( U_{(1,0)}^{\mathit{CT}}(t)
))_t$ and $( e^{-mt}u_2(U_{(0,1)}^{\mathit{CT}}(t)
))_t$ are
martingales having, respectively, $u_2(U_{(1,0)}^{\mathit{CT}}(0)
)=b/S$ and $u_2(U_{(0,1)}^{\mathit{CT}}(0))=-c/S$ as expectations.
They converge in ${\rm L}^{p}$ for every nonnegative integer $p\geq1$.
We are interested in the probability distributions of
%
\begin{equation}\label{defXY}
X:= \lim_{t\rightarrow+\infty}e^{-mt}u_2\bigl(U_{(1,0)}^{\mathit{CT}}(t) \bigr)
\quad\mbox{and}\quad  Y:= \lim_{t\rightarrow+\infty
}e^{-mt}u_2\bigl(U_{(0,1)}^{\mathit{CT}}(t) \bigr).
\end{equation}
Projecting along the second eigenline, scaling and passing to the limit in
system (\ref{dislocationVectorielle}) lead straightforwardly to the
following proposition.

\begin{Prop}
The limit random variables $X$ and $Y$ are solution of the following
(scalar) limit dislocation equations:
%
\begin{equation}
\label{dislocationProjetee}
\cases{
X\stackrel{\mathcal{L}}{=}e^{-m\tau}( [a+1]X+[b]Y), \cr
Y\stackrel{\mathcal{L}}{=}e^{-m\tau}([c]X+[d+1]Y),
}
\end{equation}
with
%
\begin{equation}\label{esperances}
\mathbb{E}(X) = \frac b S, \qquad \mathbb{E}(Y) = - \frac{c} S,
\end{equation}
where all the mentioned variables are independent.
\end{Prop}

\begin{Rem}
Janson \cite{Jan} in his Theorem 3.9 gets the same limit dislocation
equations. He obtains the unicity of the solution in $L^2$ by a fixed
point method. Hereunder in Section \ref{calculdeF}, calculating
explicitly the solution of the fixed point system (\ref
{dislocationProjetee}) together with conditions (\ref{esperances}), we
give in passing another proof of the unicity in $L^2$.
\end{Rem}


\section{Characteristic functions: fundamental differential system}

Let $\mathcal{F}$ and $\mathcal{G}$ be respectively, the characteristic
functions of $X$ and $Y$:
\[
\forall x\in\mathbb{R}\qquad
 \mathcal{F}(x)=\mathbb{E}(e^{ixX})=\int_{-\infty}^{+\infty
}e^{ixt}\,d\mu_X(t)
\]
with a similar formula for $\mathcal{G}$.
Since $X$ and $Y$ admit moments of all orders, $\mathcal{F}$ and
$\mathcal{G}$
are infinitely differentiable
on $\mathbb{R}$.

\begin{Prop}
The characteristic functions $\mathcal{F}$ and $\mathcal{G}$ are
solutions of
the differential system
%
\begin{equation}
\label{systemeFourier}
\cases{
\mathcal{F}(x)+mx\mathcal{F}'(x)=\mathcal{F}(x)^{a+1}\mathcal
{G}(x)^b,
\cr
\mathcal{G}(x)+mx\mathcal{G}'(x)=\mathcal{F}(x)^{c}\mathcal{G}(x)^{d+1},
}
\end{equation}
and satisfy the boundary conditions at the origin
%
\begin{equation}
\label{conditionsBordSystemeFourier}
\cases{
\mathcal{F}(x)=1+i\dfrac{b}Sx+O(x^2),
\cr
\mathcal{G}(x)=1-i\dfrac{c}Sx+O(x^2).
}
\end{equation}
\end{Prop}

\begin{pf}
Conditioning on $\tau$ the distribution of which is exponential with
mean~$1$, the first dislocation equation
(\ref{dislocationProjetee}) implies successively that, for any $x\in
\mathbb{R}$,
\begin{eqnarray*}
\mathcal{F}(x)&=&\mathbb{E}\bigl(
\mathbb{E}\bigl( \exp\bigl( ixe^{-m\tau}([a+1]X+[b]Y)|\tau\bigr)
\bigr)\bigr)
\\
&=&
\int_0^{+\infty}\mathcal{F}^{a+1}\bigl( xe^{-mt}\bigr)\mathcal
{G}^{b}
\bigl( xe^{-mt}\bigr)e^{-t}\,dt.
\end{eqnarray*}
After a change of variable under the integral, this functional equation
can be written
\[
\forall x\neq0\qquad
\mathcal{F}(x)=\frac{x}{m|x|^{1+1/m}}\int_0^x\mathcal
{F}^{a+1}(t)\mathcal{G}^b(t) \frac{dt}{|t|^{1-1/m}}.
\]
Differentiation of this equality and the similar one obtained from the
second dislocation equation in
(\ref{dislocationProjetee}) lead to the result.
The boundary conditions come elementarily from the computation of the
means of $X$ and $Y$ and from the existence of their second moment
(Taylor expansion of $\mathcal{F}$ and $\mathcal{G}$ at $0$).
\end{pf}

\begin{Rem}
The differential system (\ref{systemeFourier}) is singular at $0$ so
that the unicity of its solution
that satisfies the boundary condition (\ref
{conditionsBordSystemeFourier}) is not an elementary
consequence of general theorems for ordinary differential equations.
\end{Rem}


\section{Resolution of the fundamental differential system}\label{resolution}

\subsection{Change of functions: Heuristics}\label{heuristics}

Formally, without carefully checking which $m$th roots should be
considered, if the variables $x\in\mathbb{R}$
and $w\in\mathbb{C}$ are related by $x^Sw^m=1$,
the change of functions
\[
\cases{
f(w)
=x^{1/m}\mathcal{F}(x),\cr
g(w)=x^{1/m}\mathcal{G}(x)
}
\]
reduces the problems (\ref{systemeFourier}) and (\ref
{conditionsBordSystemeFourier}) to the
regular differential system
%
\begin{equation}
\label{systemeMonomial}
\cases{
f'=\dfrac{-1}{S}f^{a+1}g^b,
\cr
g'=\dfrac{-1}{S}f^cg^{d+1},
}
\end{equation}
with boundary conditions at infinity
%
\begin{equation}
\label{conditionsBordSystemeMonomial}
\cases{
f(w)=w^{-1/S}+i\dfrac bSw^{-{(1+m)}/S}+O\bigl( |w|^{-
{(1+2m)}/S}\bigr),
\vspace*{1pt}\cr
g(w)=w^{-1/S}-i\dfrac cSw^{-({1+m})/S}+O\bigl( |w|^{-
({1+2m})/S}\bigr).
}
\end{equation}


The basic fact for the resolution of (\ref{systemeMonomial}) is that
it admits
$1/g^m-1/f^m$ as first integral:
if $K$ is any complex number such that the constant function
$1/g^m-1/f^m$ equals
$1/K^m$, then $g^m$ can be straightforwardly expressed as a function of
$f$ and (\ref{systemeMonomial})
implies that $f$ is solution of the ordinary differential equation
%
\begin{equation}
\label{equaDiffAbelienne}
f'\times\frac{( 1+( f/K)^m
)^{b/m}}{(  f/K)^{S+1}}
=-\frac{K^{S+1}}{S}
\end{equation}
with boundary conditions coming from (\ref{conditionsBordSystemeMonomial}).

This leads to consider a primitive of the function $z\mapsto
(1+z^m)^{b/m}/z^{S+1}$ in the
complex field.


\subsection{Abelian integral $I$ and its inverse $J$}

For all integers $m$, $S$ and $b$ that satisfy $S\geq5$, $S/2< m< S$,
$1\leq b<S/2$, we denote by
$I=I_{m,S,b}$ the function
\[
I(z)=\int_{[z,z\infty)}(1+u^m)^{ b/m}\frac
{du}{u^{S+1}}
=\frac1{z^S}\int_1^{+\infty}[ 1+(tz)^m]
^{ b/m}\frac{dt}{t^{S+1}},
\]
where $[z,z\infty)$ denotes the ray $\{ tz, t\geq1\}$ and where the
power $1/m$ is used for the
principal determination of the $m$th root.
The function $I$ is an Abelian integral on the curve $x^m-y^m=1$
(which is isomorphic to the famous Fermat curve $x^m+y^m=1$ by a
straightforward linear change
of variables),
defined on the open set
\[
\mathcal{O}_m=\mathbb{C}\bigm\backslash\bigcup_{p\in\{ 0,\dots,m-1\}
}\mathbb{R}_{\geq
0}e^{({i\pi}/m)(1+2p)}.
\]
Note that the integral is convergent because $S-b+1\geq3$.
Let $\mathcal{S}_m$ be the open sector of the complex plane defined by
\[
\mathcal{S}_m=\biggl\{ z\in\mathbb{C}\setminus\{ 0\},\ -\frac\pi m<\arg
(z)<\frac\pi
m\biggr\}.
\]
The open set $\mathcal{O}_m$ is the union of the images of $\mathcal{S}_m$
under all rotations of angles $2k\pi/m$
around the origin, $k\in\mathbb{Z}$.

In the following, the notation ${b/m\choose n}$ denotes the ordinary
binomial coefficient, generalized
for rational (or even complex) values of $b/m$ by Euler's Gamma function.
As everywhere else in the paper, the positive integer $a$ is $a=S-b$.

\begin{Prop}[(Properties of $I$)]\label{developI}
\begin{longlist}[(iii)]
\item[(i)]
$I$ is holomorphic on $\mathcal{O}_m$ and for any $z\in\mathcal{O}_m$,
%
\begin{equation}
\label{primitiveI}
I'(z)=-\frac{( 1+z^m) ^{ b/m}}{z^{S+1}}.
\end{equation}

\item[(ii)]
For any $m$th root of unity $\omega$ and for any $z\in\mathcal{O}_m$,
%
\begin{equation}
\label{rotationsI}
I(\omega z)=\omega^{-S}I(z).
\end{equation}

\item[(iii)]
The function $I$ admits a power series expansion in the neighborhood of
infinity in any connected component of $\mathcal{O}_m$.
On $\mathcal{S}_m$, this expansion is given by the formula
%
\begin{equation}
\label{DSEinfinideI}
\qquad I(z)=\sum_{n\geq0}\frac1{a+mn}\pmatrix{ b/m\cr n }z^{-a-mn}
=\frac1{az^a}+\frac b{m(a+m)}\frac1{z^{a+m}}+\cdots,
\end{equation}
valid for any $z\in\mathcal{S}_m$, $|z|\geq1$.

\item[(iv)]
The function $I$ admits a Laurent series expansion in the neighborhood
of the origin in any connected component of $\mathcal{O}_m$.
On $\mathcal{S}_m$, this expansion is given by the formula
%
\begin{equation}
\label{DSL0deI}
I(z)=\frac1{Sz^S}+\frac b{m(S-m)}\frac1{z^{S-m}}
+C_0
-\sum_{n\geq2}\pmatrix{ b/m\cr n }\frac{z^{mn-S}}{mn-S},
\end{equation}
where $C_0$ is the constant
%
\begin{equation}
\label{serieC_0}
C_0=\sum_{n\geq0}\pmatrix{ b/m\cr n }\biggl( \frac1{a+mn}+\frac
1{mn-S}\biggr) .
\end{equation}
Formula (\ref{DSL0deI}) is valid for any $z\in\mathcal{S}_m$,
$|z|\leq1$.

\item[(v)]
$C_0<0$.
\end{longlist}
\end{Prop}

\begin{pf}
(i) and (ii) are direct consequences of the definition of $I$.
Expansion~(\ref{DSEinfinideI}) and its validity for $z\in\mathcal
{S}_m$, $|z|>1$
comes directly from the power series expansion of $\zeta\mapsto
(1+\zeta)^{b/m}$ in the definition of $I$.
Its validity for $|z|=1$ is given by the convergence of the series at
such a $z$ and application
of Abel's theorem,\footnote{We refer to the following theorem of Abel: if a series $\sum_na_n$ is
convergent, then the power
series $\sum_na_nz^n$ converges uniformly on the segment $[0,1]$.} proving (iii).
To prove expansion (\ref{DSL0deI}), notice first that $I$ is
holomorphic on the simply connected
domain $\mathcal{S}_m$ and $I'(z)$ tends to $0$ as $z$ tends to infinity,
so that integration on the
ray $[z,z\infty)$ is equivalent to integration on $[z,1]$ followed
by integration on $[1,+\infty)$.
Thus,
\[
I(z)=I(1)+\int_{[z,1]}(1+u^m)^{b/m}\frac{du}{u^{S+1}}.
\]
Power series expansion of $u\mapsto(1+u)^{b/m}$ under this last
integral leads then to~(\ref{DSL0deI}).
The proof of (iv) is again made complete by application of Abel's theorem.
Note that, since $S$ is not a multiple of $m$ because of our
assumptions on the parameters, the
denominators in Formula (\ref{DSL0deI}) do not vanish.
Finally, if $\alpha_n$ denotes the general term of the series (\ref
{serieC_0}), a straightforward
computation shows that
\[
\alpha_0+\alpha_1=
\frac{(S-a)(m^2+aS)(S-a-m)}{amS(a+m)(S-m)}<0,
\]
the last inequality coming from $S-a-m<S-S/2-S/2=0$ and from the other
hypotheses on the parameters.
Furthermore, $\alpha_{2n}+\alpha_{2n+1}<0$ for any $n\geq1$, which
concludes the proof.
[Hint: compute $\alpha_{2n}+\alpha_{2n+1}$, factorize $
{{b/m} \choose {2n}}$ by ${{b/m} \choose {2n+1}}$,
use the fact $(2n+1)/(2n-b/m)>1$, notice that ${{b/m}\choose {2n+1}}>0$ because
$0<b/m=(S-a)/m<S/2m<1$.]
\end{pf}

Let $\mathbb{H}$ denote Poincar\'e half-plane:
\[
\mathbb{H}=\{ z\in\mathbb{C}, \Im(z)>0\} \quad  \mbox{and}
\quad
\overline{\mathbb{H}} = \{ \overline z, z \in\mathbb{H}\} .
\]

\begin{figure}

\includegraphics{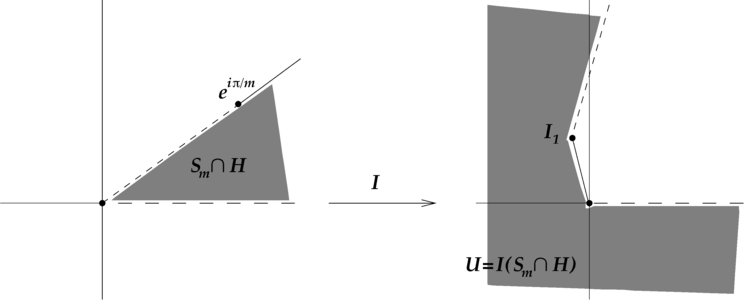}

\caption{Domain $\mathcal{S}_m\bigcap\mathbb{H}$ and its image by
$I$.}\label{fig2}
\end{figure}

\begin{Prop}
\label{Iconforme}
The analytic function $I\dvtx \mathcal{S}_m\bigcap\mathbb{H}\to\mathbb
{C}$ is a conformal
mapping onto the open subset
\[
\mathcal{U}=\biggl\{ z, -\frac{a\pi}{m}<\arg(z)<0\biggr\}
\cup\biggl( I_1 +\biggl\{ z,\ -\frac{S\pi}m<\arg(z)<-\frac{a\pi}{m}\biggr\}
\biggr)
\]
(see Figure \ref{fig2}), where
%
\begin{equation}
\label{I1}
I_1:=
\frac1m B\biggl(\frac am,\frac dm\biggr)e^{-{(ia\pi)/m}},
\end{equation}
and where $B$ denotes Euler's Beta function $B(x,y)=\Gamma(x)\Gamma
(y)/\Gamma(x+y)$.
\end{Prop}

\begin{pf}
Let $\zeta_m = \exp(i\pi/m)$.
We show hereunder that the restriction of $I$ to the sector $\mathcal
{S}_m\cap Cl(\mathbb{H})$ (where $Cl(\mathbb{H})$ denotes the
topological closure of
$\mathbb{H}$) admits a continuous
continuation to the ray $\{ t \zeta_m, t>0\}$
and that this continuation maps homeomorphically the
boundary of the sector $\mathcal{S}_m\cap\mathbb{H}$ onto the
boundary of $\mathcal{U}$.
The result is then a consequence of elementary geometrical conformal theory
(see, for example, Saks and Zygmund \cite{SaksZygmund}).

Let $h\in\mathbb{H}$, $r>0$, $t>1$ and $z=r(1-h)\zeta_m$.
When $h$ tends to $0$, then $1+(tz)^m=1-r^mt^m+mr^mt^mh+O(h^2)$ so that
the value
of $m$th root principal determination of $1+(tz)^m$ according to the
sign of $1-(rt)^m$ leads to the
respective limits in terms of Beta incomplete functions:

\begin{itemize}[$\bullet$]
\item[$\bullet$] if $r\geq1$, then
%
\begin{equation}
\label{IBordInf}
\lim_{z\to r\zeta_m, z\in\mathcal{S}_m}I(z)=
\frac1m\zeta_m^{-a}\int_0^{1/r^m}(1-u)^{b/m}u^{c/m}\,du;
\end{equation}

\item[$\bullet$] if $r\leq1$, then
%
\begin{equation}
\label{IBordSup}
\lim_{z\to r\zeta_m, z\in\mathcal{S}_m}I(z)=
I_1+\frac1m\zeta_m^{-S}\int_1^{1/r^m}(u-1)^{b/m}u^{c/m}\,du.
\end{equation}
\end{itemize}
The complex number $I_1$ is simply
\[
I_1=\lim_{z\to\zeta_m, z\in\mathcal{S}_m}I(z);
\]
formula (\ref{I1}) is a consequence of the integral representation of
Euler Beta function
$B(\alpha, \beta)=\int_0^1(1-u)^{\alpha-1}u^{\beta-1}\,du$.
The monotonicity of real integrals (\ref{IBordInf}) and (\ref
{IBordSup}) with respect to $r$ show that
the continuous continuation of $I$ defined by these formulae maps
decreasingly the ray $]0,+\infty[$
onto itself and respectively, the ray $]0,\zeta_m]$ onto the ray $\{
I_1+ t\zeta_m^{-S}, t\geq0 \}$
and the ray $[\zeta_m,\zeta_m\infty)$ onto $[I_1,0[$.
\end{pf}

\begin{Rem}
By computation in the realm of hypergeometric functions, one shows that
the numbers $C_0$ defined
by (\ref{serieC_0}) and $I_1$ defined by (\ref{I1}) are related by
\[
C_0=-\frac{\sin\pi(1+b/m)}{\sin\pi(1+S/m)}|I_1|
=-\frac1m\frac{\sin\pi(1+b/m)}{\sin\pi(1+S/m)}B\biggl( \frac
{S-b}m,\frac{m+b}m\biggr) .
\]
\end{Rem}

\begin{Def}
Let $J=J_{m,S,b}\dvtx \mathbb{C} \setminus]-\infty,0]\to\mathcal{S}_m$
the only
continuous function defined by:
\begin{itemize}
\item[$\bullet$] $\forall z\in\overline{\mathbb{H}}$, $J(z)=I^{-1}(z)$ in the
sense of Proposition \ref{Iconforme}
($\overline{\mathbb{H}}$ is an open subset of $\mathcal{U}$
so that this
functional inverse exists);

\item[$\bullet$] $\forall z\in\mathbb{H}$, $J(z)=\overline{J(\overline z)}$
(complex conjugacy).
\end{itemize}
\end{Def}

\begin{figure}[t]

\includegraphics{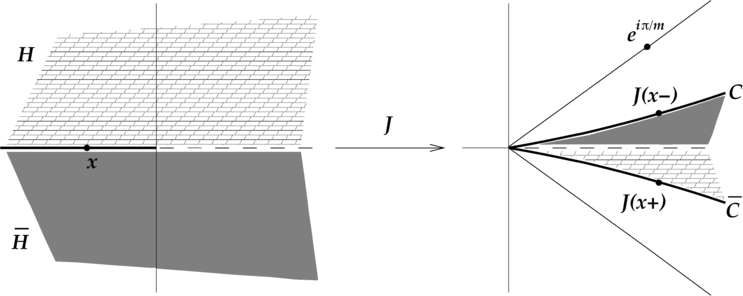}

\caption{Action of $J$ on the slit plane $\mathbb{C}\setminus\mathbb
{R}_{\leq0}$.}\label{fig3}
\end{figure}

The properties of $I$ shown in Propositions \ref{developI} and \ref
{Iconforme} imply that $J$ is
a conformal mapping between $\mathbb{C} \setminus]-\infty,0]$ and an open
subset of $\mathcal{S}_m$
(use Schwarz reflection principle), that maps
$\mathbb{H}$ into $\mathcal{S}_m\cap\overline{\mathbb{H}}$ and
$\overline{\mathbb{H}}$
into $\mathcal{S}_m\cap{\mathbb{H}}$.
If $\mathcal{C}$ denotes the inverse of the negative real axis by the
restriction of $I$ to $\mathcal{S}_m\cap\mathbb{H}$,
then the boundary of the image of $J$ is $\mathcal{C}\cup\overline
{\mathcal{C}}\cup\{ 0\}$ (see Figure \ref{fig3}).
Furthermore, the restriction of $J$ to the positive real half-line is
the inverse of $I$'s and $J$ is the
unique analytic expansion of $(I_{|]0,+\infty[})^{-1}$ to the slit plane.
Naturally, the formula $J(\overline z)=\overline{J(z)}$ is valid when
$z$ is any nonnegative complex
number.

\begin{Prop}
\label{imageJ}
For any negative real number $x$, both limits
\[
\lim_{z\to x, z\in\mathbb{H}}J(z)
 \quad \mbox{and}\quad
\lim_{z\to x, z\in\overline{\mathbb{H}}}J(z)
\]
exist, are nonreal and conjugate (thus, different).
\end{Prop}

\begin{pf}
Direct consequence of the preceding properties of $J$ and Proposition~\ref{Iconforme} (see Figure \ref{fig3}).
\end{pf}

We adopt the following notation:
%
\begin{equation}
\label{Jlimite}
\forall x<0\qquad
\cases{
J(x-)=\displaystyle\lim_{z\to x, z\in\overline{\mathbb{H}}}J(z)\in\mathcal
{S}_m\cap\mathbb{H},
\cr
J(x+)=\displaystyle\lim_{z\to x, z\in\mathbb{H}}J(z)\in\mathcal{S}_m\cap
\overline{\mathbb{H}}.
}
\end{equation}

\begin{Prop}
\label{DSPJ}
The function $J$ admits, as $z$ tends to infinity in the slit plane
$\mathbb{C}\setminus\mathbb{R}_-$, an
asymptotic Puiseux series expansion at any order in the scale
\[
\biggl(\frac1z\biggr)^{1/S+p\sigma+q},\qquad (p,q)\in\mathbb{N}^2,
\]
where all fractional powers denote principal determination.
The beginning of this asymptotic expansion is
\begin{eqnarray}
\label{DSPinfiniDeJ}
J(z)&=&\biggl( \frac1{Sz}\biggr)^{1/S}
+\frac{b}{m(S-m)}\biggl( \frac1{Sz}\biggr)^{{(m+1)}/{S}}\nonumber
\\[-8pt]\\[-8pt]
&&{}+C_0\biggl( \frac1{Sz}\biggr)^{{(S+1)}/{S}}
+o\biggl(\frac{1}{z}\biggr)^{{(S+1)}/{S}}.\nonumber
\end{eqnarray}
\end{Prop}

\begin{pf}
Expansion (\ref{DSL0deI}) leads to (\ref{DSPinfiniDeJ}) using the
reversion formula $J\circ I=\break \operatorname{Id}$.
\end{pf}


\subsection{Computation of characteristic functions}
\label{calculdeF}

This section gives an explicit closed form of characteristic functions
$\mathcal{F}$ and $\mathcal{G}$ for
the elementary continuous-time urn processes $X$ and $Y$ [defined in
(\ref{defXY})] associated with the replacement matrix
$R=\left(
{a\atop c}\enskip {b\atop d
}
\right)$, in terms of the just defined functions $J$.
Remember: the urn is supposed to be large and nontriangular so that
$b>0$ and $c> 0$.
Let $\kappa$ be the positive number defined by
%
\begin{equation}
\label{kappa}
\kappa
=  \sqrt[m\,]{\frac{S}{m(S-m)}}.
\end{equation}

\begin{Th}
\label{characteristicFunctions}
The characteristic functions $\mathcal{F}$ and $\mathcal{G}$ are the unique
solutions of the differential system
(\ref{systemeFourier}) that satisfy boundary conditions (\ref
{conditionsBordSystemeFourier}).
They are given by the formulae
%
\begin{equation}
\label{formulesFG}
\hspace*{15pt} \forall x>0\hspace*{8pt}
\cases{
\mathcal{F} (x)=\kappa e^{-{i\pi/(2m)}}x^{-1/m}
J_{m,S,b}\biggl( C_0+\dfrac{\kappa^Se^{-{i\pi
S/(2m)}}}{S}x^{-S/m}\biggr),
 \cr
\mathcal{G} (x)=\kappa e^{{i\pi/(2m)}}x^{-1/m}
J_{m,S,c}\biggl( C_0+\dfrac{\kappa^Se^{{i\pi
S/(2m)}}}{S}x^{-S/m}\biggr)
}\hspace*{-12pt}
\end{equation}
and
%
\begin{equation}
\label{moitie}
\forall x\in\mathbb{R}\qquad \mathcal{F}(-x)=\overline{\mathcal{F}(x)},
\qquad \mathcal{G}(-x)=\overline{\mathcal{G}(x)}.
\end{equation}
\end{Th}

\begin{pf}
(1)
We first solve (\ref{systemeFourier}) on $\mathbb{R}_{>0}$.
Let $F$ and $G$ be solutions of (\ref{systemeFourier}) that satisfy
(\ref{conditionsBordSystemeFourier}).
Lets do the change of variable $x\in\mathbb{R}_{>0}\to w=x^{-S/m}\in
\mathbb{R}_{>0}$ and the change of functions
\[
f(w)=w^{-1/S}F( w^{-m/S})
\quad \mbox{and} \quad
g(w)=w^{-1/S}G( w^{-m/S})
\]
that is straightforwardly reversed by the formula
$F(x)=x^{-1/m}f(x^{-S/m})$ and a similar one for
$G$ and $g$.
Then $f$ and $g$ are solutions of (\ref{systemeMonomial}) on $\mathbb
{R}_{>0}$ and satisfy boundary
conditions (\ref{conditionsBordSystemeMonomial}) at $+\infty$.
In particular, since (\ref{systemeMonomial}) is a nonsingular
differential system, Cauchy--Lipschitz
theorem guarantees that if $(f,g)$ is any solution, then $f$
({resp.}, $g$) is identically zero or does
not vanish.
This implies that $f$ and $g$ do not vanish on $\mathbb{R}_{>0}$.
Because of the balance conditions $a+b=c+d$, differentiation of
$1/g^m-1/f^m$ leads to the fact that
this function is constant  on $\mathbb{R}_{>0}$ (first integral).
Furthermore, boundary conditions at $+\infty$ (\ref
{conditionsBordSystemeMonomial}) imply that this
constant value is $i\frac mS(b+c)$.
If $K$ denotes the complex number
\[
K=\kappa\exp\biggl(-\frac{i\pi}{2m}\biggr)
\]
[$\kappa>0$ has been defined by formula (\ref{kappa})], this shows that
%
\begin{equation}
\label{integralePremiere}
\forall w>0\qquad \frac1{g^m(w)}-\frac1{f^m(w)}=\frac1{K^m}.
\end{equation}
Since $f/g$ is continuous on $\mathbb{R}_{>0}$, does not vanish and
tends to
$1$ at $+\infty$
(\ref{conditionsBordSystemeMonomial}), relation $(f/g)^m=1+(f/K)^m$
implies that, on $\mathbb{R}_{>0}$,
%
\begin{equation}
\label{gFonctionDef}
g=\frac{f}{(1+( {f/K}) ^m )^{1/m}}
\end{equation}
(principal determination of the $m$th root).
Reporting in the first equation of (\ref{systemeMonomial}) shows that
$f$ is necessarily a solution
of equation (\ref{equaDiffAbelienne}) on $\mathbb{R}_{>0}$.
Boundary conditions (\ref{conditionsBordSystemeMonomial}) imply that,
when $w$ tends to $+\infty$, $\frac1Kf(w)\sim\frac1\kappa e^{i\pi
/2m}w^{-1/S}\in\mathcal{S}_m$, so that
equation (\ref{equaDiffAbelienne}) can be written
\[
\frac d{dw}I_{m,S,b}\biggl(\frac{f(w)}K\biggr)=\frac{K^S}S
\]
in a neighborhood of $+\infty$.
Integration of this equation shows that
\[
I_{m,S,b}\circ\biggl( \frac fK\biggr) (w)=\frac{K^S}Sw+C_1
\]
in a neighborhood
of $w=+\infty$, for a suitable complex constant $C_1$.
The determination of $C_1$ is made by means of local expansions:
since $f$ tends to $0$ at $+\infty$, using (\ref{DSL0deI}) and
previous equality leads to
\[
C_1+\frac{K^S}Sw
=\frac{K^S}{Sf(w)^S}
+\frac{b}{m(S-m)}\frac{K^{S-m}}{f(w)^{S-m}}
+C_0
+o(1),
\]
when $w$ tends to $+\infty$, so that boundary conditions (\ref
{conditionsBordSystemeMonomial})
lead to the equality $C_1=C_0$.
Note that this computation makes use of the big-O in (\ref
{conditionsBordSystemeMonomial}), of the
assumption $1-2m/S<0$ (large urn) and of the relation $S-m=b+c$.
Thus, necessarily,
%
\begin{equation}
\label{formulef}
f(w)=KJ_{m,S,b}\biggl( C_0+\frac{K^S}Sw\biggr)
\end{equation}
for any $w$ in a neighborhood of $+\infty$.
The function $w\to KJ_{m,S,b}(C_0+K^Sw/S)$ is well defined on $\mathbb
{R}_{>0}$ because
$C_0<0$ [Proposition \ref{developI}(5)] and $-\pi
<\operatorname{arg}
(K^S)<-\pi/2$, so that it is the only maximal
solution on $\mathbb{R}_{>0}$ of equation (\ref{equaDiffAbelienne}) that
satisfies the first equation of
(\ref{conditionsBordSystemeMonomial}).
This shows finally that
\[
\forall x>0\qquad
F (x)= Kx^{-1/m}
J_{m,S,b}\biggl( C_0+\frac{K^S}Sx^{-S/m}\biggr).
\]
Since $-K^m=\overline K^m$, the same arguments show that, for any $w>0$,
\[
g(w)=\overline KJ_{m,S,c}\biggl( C_0+\frac{\overline K^S}{S}w\biggr) ,
\]
which shows completely formula (\ref{formulesFG}).

(2)
The resolution on $\mathbb{R}_{<0}$ is made the same way.
To this effect, lets do the new change of variable
$x\in\mathbb{R}_{<0}\to w=|x|^{-S/m}e^{i\pi S/m}\in\mathbb
{R}_{>0}e^{i\pi S/m}$.
Lets do as well the change of functions
\[
f(w)=e^{-i\pi/m}|w|^{-1/S}F( -|w|^{-m/S})
\]
 and
\[
g(w)=e^{-i\pi/m}|w|^{-1/S}G( -|w|^{-m/S}).
\]
These changes of variable and functions are reversed by the formulae
$x=-|w|^{-m/S}$ and
$F(x)=e^{i\pi/m}|x|^{-1/m}f(e^{i\pi S/m}|x|^{-S/m})$ with a similar
formula for $G$ and $g$.
Functions $f$ and $g$ are still solutions of (\ref{systemeMonomial})
but boundary conditions become,
as $|w|$ tends to infinity,
%
\begin{equation}
\label{conditionsBordSystemeMonomialNegatif}
\cases{
f(w)=e^{-i{(\pi/m)}}|w|^{-1/S}
\biggl( 1-i\dfrac bS|w|^{-{m}/S}+O( |w|^{-{2m}/S}
)\biggr),
\cr
g(w)=e^{-i({\pi}/{m})}|w|^{-1/S}
\biggl( 1+i\dfrac cS|w|^{-{m}/S}+O( |w|^{-{2m}/S}
)\biggr) .
}
\end{equation}
This implies that First integral (\ref{integralePremiere}) is still
valid (same $K$) and, since $f$ and
$g$ are still equivalent at infinity, relation (\ref{gFonctionDef}) is
satisfied.
Boundary conditions (\ref{conditionsBordSystemeMonomialNegatif}) imply that,
when $w$ tends to $+\infty$, $\frac1Kf(w)\sim\frac1\kappa
|w|^{-1/S}e^{-i\pi/2m}\in\mathcal{S}_m$.
Consequently, the same arguments as before show that formula (\ref
{formulef}) remains valid
(note that $C_0+wK^S/S\in\mathbb{H}$ so that this formula is well defined
for any $w$).
This shows that
\[
\forall x<0\qquad
F (x)= Ke^{{i\pi}/{m}}|x|^{-1/m}
J_{m,S,b}\biggl( C_0+\frac{K^S}Se^{i\pi({S}/{m})}|x|^{-
S/m}\biggr).
\]
Since $Ke^{i\pi/m}=\overline K$, one gets finally $F(-x)=\overline
{F(x)}$ for any real number $x$.
The proof of the whole theorem is made complete by the same arguments
for $G$.
\end{pf}

\begin{Rem}
Formula (\ref{moitie}) on characteristic functions comes directly from
the fact that $X$ and $Y$ are
real-valued random variables.
\end{Rem}

We want to know more about the analyticity properties of $\mathcal{F}$ and
$\mathcal{G}$ around $0$.
Let $\varphi=\varphi_{m,S,b}$ be the function defined by the formula
%
\begin{equation}
\label{defPhi}
\varphi(z)=\kappa z^{-1/m}J_{m,S,b}\biggl( C_0+\frac{\kappa
^S}{S}( z^{-1/m}) ^S\biggr),
\end{equation}
where the power $1/m$ denotes the principal determination of the $m$th root.
Note that $\kappa$ and $C_0$, respectively, defined by formulas (\ref
{kappa}) and (\ref{serieC_0})
are functions of $m$, $S$ and $b$ too.
If $\rho$ denotes the positive number
\[
\rho= \biggl( \frac{S|C_0|}{\kappa^S}\biggr) ^{-m/S}
=\frac{S^{1-S/m}|C_0|^{-m/S}}{m(S-m)},
\]
it follows from the properties of $J$ that $\varphi$ is defined and
holomorphic on the open set
\[
\mathcal{V}=\mathbb{C}\setminus\{ (-\infty,0]\cup[\rho
,+\infty)
\}.
\]
Furthermore, the characteristic functions $\mathcal{F}$ and $\mathcal
{G}$ are
restrictions of $\varphi$ functions
on the imaginary axis:
for any $x\in\mathbb{R}$,
\[
\mathcal{F}(x)=\varphi_{m,S,b}(ix)
\quad \mbox{and}\quad
\mathcal{G}(x)=\varphi_{m,S,c}(-ix).
\]
Note that $\kappa$ is a function of $(m,S)$ so that the same
$\kappa$ appears in both functions $\varphi_{m,S,b}$ and $\varphi_{m,S,c}$
(respective constants $C_0$ and $\rho$ are however different).

\begin{Prop}\label{analyticPhi}
The function $\varphi$, holomorphic on $\mathcal{V}$, cannot be
analytically extended on a larger subset
of $\mathbb{C}$.
However, setting $\varphi(0)=1$ defines a
continuously differentiable extension of $\varphi$ on $\mathcal
{V}\cup\{
0\}$.
\end{Prop}

\begin{pf}
The half-line $[\rho,+\infty)$ is the locus of complex $z$ such that
$C_0+\frac{\kappa^S}{S}( z^{-1/m}) ^S$ is a real
nonpositive number
(remember that $m<S<2m$).
Since the principal determination of the $m$th root is well defined and
nonzero in a neighbourhood of
this half-line, Proposition \ref{imageJ} implies that $\varphi$
cannot be continuously extended at any
point of $[\rho,+\infty)$.

If $x$ is a negative number, definition of the principal determination
of the $m$th root leads to the
existence of both limits
\[
\cases{\displaystyle
\lim_{z\to x, z\in\mathbb{H}}\varphi(z)=
\kappa e^{-i{(\pi/m)}}|x|^{-1/m}J\biggl( C_0+\dfrac{\kappa
^S}{S}e^{-i{(\pi S/m)}}|x|^{- S/m}\biggr)
:=\varphi(x+),
\cr
\displaystyle\lim_{z\to x, z\in\overline{\mathbb{H}}}\varphi(z):=\varphi
(x-)=\overline{\varphi(x+)}.
}
\]
Since
the image of $J$ is included in $\mathcal{S}_m$, the limit $\varphi(x+)$
belongs to the
open sector $e^{-i{(\pi/m)}}\mathcal{S}_m$ which contains no real
number, so that
$\varphi(x+)\neq\varphi(x-)$.
This shows that $\varphi$ cannot be continuously extended at any point
of $\mathbb{R}_{<0}$.

When $z$ tends to $0$ in the slit plane $\mathbb{C}\setminus\mathbb{R}_{<0}$,
Proposition \ref{DSPJ} shows
that $\varphi(z)$ tends to $1$.
One step more, computing the derivative of $\varphi$ in terms of $J$
using the algebraic expression
of $I'$ (\ref{primitiveI}) implies, with expansion (\ref
{DSPinfiniDeJ}), that
\[
\lim_{z\to0, z\in\mathbb{C}\setminus\mathbb{R}_{\leq0}}\varphi
'(z)=\frac bS.
\]
\upqed\end{pf}

\begin{Cor}\label{rayonnul}
The exponential moment generating series
\[
\sum\frac{\mathbb{E}(X^p)}{p!}T^p
  \quad\mbox{and}\quad
\sum\frac{\mathbb{E}(Y^p)}{p!}T^p
\]
have a radius of convergence equal to $0$.
\end{Cor}

\begin{pf}
These series are the Taylor series of $\varphi_{m,S,b}$ and $\varphi
_{m,S,c}$ at $0$.
If these radii were positive, these functions could be analytically
extended to a neighborhood of the
origin.
\end{pf}

\begin{Rem}
The singularity of $\varphi$ at the origin is thus not due to
ramification but to a divergent Taylor series
phenomenon.
Indeed, the apparent ramification coming from the $m$th root at the
origin in formula (\ref{defPhi})
is compensated by both Puiseux expansion (\ref{DSPinfiniDeJ}) and the
$S$th power of the $m$th root appearing in the argument of $J$ in
formula (\ref{defPhi}).
\end{Rem}


\section{Density of $W^{\mathit{CT}}$}
\label{densiteWCT}

Notice, with the notation (\ref{Jlimite}) that
%
\begin{equation}
\label{equiv}
\mathcal{F}(x)\mathop{\sim}\limits_{x\to+\infty}\kappa
J(C_0-)x^{-1/m},
\end{equation}
where the nonreal complex number $J(C_0-)$ is different from $0$ (see Figure 3).

A first consequence is that $\mathcal{F}(x)$ tends to $0$ when $x$
tends to
$+\infty$.
Hence, the probability distribution function of $W^{\mathit{CT}}$ is continuous
so that the law of $W^{\mathit{CT}}$ has no point mass.

A second consequence is that $\mathcal{F}$ is not in ${\rm L}^1$ so
that $W^{\mathit{CT}}$
distribution cannot be obtained by classical Fourier inversion.
Nevertheless, we obtain in Section \ref{inversionFourier} an
expression of this density using the derivative of the characteristic
function $\mathcal{F}$.
Before, we need firstly to ensure that the support of $W^{\mathit{CT}}$ is the
whole real line ${\mathbb{R}}$ which is proven
in Section \ref{support} and secondly to ensure that $W^{\mathit{CT}}$ admits a
density which is proven in
Section \ref{densityconnection} using the martingale connection (\ref
{martingaleconnection}).
As usually, this kind of connection induces a smoothing phenomenon
between $W^{\mathit{DT}}$ and $W^{\mathit{CT}}$, allowing us to prove that $W^{\mathit{CT}}$ has a
density, whatever $W^{\mathit{DT}}$ distribution is.


\subsection{Support of $W^{\mathit{CT}}$}\label{support}

\begin{Prop}
\label{supportWCT}
The support of $W^{\mathit{CT}}$ is ${\mathbb{R}}$.
\end{Prop}

\begin{pf}
As in (\ref{defXY}), let $X$ denote the random variable $W^{\mathit{CT}}$
starting from one red ball.
Because of the branching property (see beginning of Section \ref
{dislocationvect}), it suffices to prove
that the support of $X$ is the whole real line $\mathbb{R}$.
General results on infinite divisibility (see, for instance, Steutel and van Harn  \cite{Steutel}, page 186) ensure that the support of an
infinitely divisible random variable having a continuous probability
distribution function is either a half-line or ${\mathbb{R}}$.
Suppose that the support of $X$ is $[\alpha, +\infty[$ for a given
real number $\alpha$.
Then denoting $X$ distribution by $\mu_X$,
\[
\mathbb{E}(e^{-sX}) =
\int_{\alpha}^{+\infty} e^{-st} \,d\mu_X(t) =
e^{-s\alpha} \int_{\alpha}^{+\infty} e^{-s(t-\alpha)} \,d\mu_X(t)
\]
exists for every real number $s\geq0$.
Hence, the function $L\dvtx s \rightarrow\mathbb{E}(e^{-sX})$ is analytic
on the
half-plane $\{ \Re z >0\}$, continuous
on the boundary of this half-plane and $\lim
_{t\rightarrow\pm\infty} \mathbb{E}(e^{itX}) = 0$.
By unicity of the analytic continuation, necessarily:
\[
L(s)
= \varphi(-s)\qquad \forall s, \Re(s) \geq0,
\]
where $\varphi$ has been introduced in (\ref{defPhi}).
But it has been proven in Proposition \ref{analyticPhi} that $\varphi
$ cannot be analytically extended on
the half-plane $\{ \Re z <0\}$.
There is a contradiction: the support of $X$ cannot be a half-line
$[\alpha, +\infty[$.

In the same way, if we suppose that the support of $W^{\mathit{CT}}$ is $]-
\infty, \beta]$ for a given real number $\beta$, we are led to a
contradiction, because $\varphi$ cannot be analytically extended on
the whole
half-plane $\{ \Re z >0\}$ (Proposition \ref{analyticPhi}).
\end{pf}


\subsection{Connection between the distribution of $W^{\mathit{DT}}$ and the
density of $W^{\mathit{CT}}$}\label{densityconnection}


\begin{Prop}
\label{density1}
Let $\mu$ be the distribution of $W^{\mathit{DT}}$ (it is a probability measure
on ${\mathbb{R}}$).

(1)
$W^{\mathit{CT}}$ admits a density $p$ on ${\mathbb{R}}$ given by
\[
\cases{
\forall w>0\qquad
p(w)=\dfrac{1}{\sigma}\dfrac{1}{\Gamma( {1/S})} w^{- 1+ \frac 1m }
\displaystyle\int_{]0,+\infty[}v^{-{1}/m} e^{-(  w/v)
^{{1}/{\sigma}}} d\mu(v),
\cr
\forall w<0\qquad
p(w)=\dfrac{1}{\sigma}\dfrac{1}{\Gamma( {1}/{S})} |w|^{- 1 +\frac 1m}
\displaystyle\int_{]-\infty,0[}|v|^{-{1}/m} e^{-( w/v)
^{{1}/{\sigma}}} \,d\mu(v).
}
\]


(2)
The density $p$ is infinitely differentiable and increasing on
${\mathbb{R}}
_{<0}$, infinitely differentiable and
decreasing on ${\mathbb{R}}_{>0}$;
it is not continuous at $0$:
$\lim_{w\rightarrow0,w\neq0} p(w)= +\infty$.
In particular, the distribution is unimodal, the mode is $0$.
\end{Prop}

\begin{pf}
(1)
To exhibit a density, let us take any real-valued bounded continuous
function $h$ defined on $\mathbb{R}$ and,
thanks to the martingale connection (\ref{martingaleconnection}), compute
\[
\mathbb{E}(h(W^{\mathit{CT}})) = \int_{{\mathbb{R}}}\int_0^{+\infty} h(uv)
g(u) \,du \,d\mu(v),
\]
where $g$ is the density of $\xi^{\sigma}$.
After the change of variable $w=uv$, we get
\begin{eqnarray*}
\mathbb{E}(h(W^{\mathit{CT}})) &=&
\int_{]-\infty, 0[}\frac{d\mu(v)}{|v|} \int_{-\infty
}^0 h(w) g\biggl( \frac{w}v\biggr)\, dw
\\
&&{}
+ \mu(\{0\})h(0)
+ \int_{]0, +\infty[}\frac{d\mu(v)}{v} \int_0^{+\infty} h(w)
g\biggl( \frac{w}v\biggr)\, dw.
\end{eqnarray*}
Recall that $W^{\mathit{CT}}$ has no point mass (see Section \ref{densiteWCT},
introductory paragraph),
so we get that $W^{\mathit{CT}}$ admits a density given by
%
\begin{eqnarray}
\label{defp}
\qquad p(w) = {\bf1}_{{\mathbb{R}}_{<0}} (w) \int_{]-\infty, 0[} g\biggl(
\frac
{w}v\biggr) \frac{d\mu(v)}{|v|}
+ {\bf1}_{{\mathbb{R}}_{>0}} (w) \int_{]0,+\infty[} g\biggl( \frac
{w}v
\biggr)\frac{d\mu(v)}v .
\end{eqnarray}
The only point to verify is that the integrals in formula (\ref{defp})
are well defined.
The density $g$ is explicit.
To simplify the notation, we consider the case when we start from one
ball ($u = 1$).
In this case,
%
\begin{equation}\label{densitedeG}
g(x) = \frac{1}{\sigma}\frac{1}{\Gamma({1}/{S})} x^{
- 1+ \frac 1m}e^{-x^{ {1}/{\sigma}}} {\bf1}_{x>0},
\end{equation}
so that, for any nonzero $w$,
\[
\frac{1}{|v|}g\biggl( \frac{w}v\biggr)
= C |w|^{- 1+\frac{1}{m}} |v|^{-{1}/m} e^{-|w|^{{1}/{\sigma
}}|v|^{-{1}/{\sigma}}}
\]
is bounded as a function of $v$.

(2) Let us prove that $\lim_{w\rightarrow0^+}
p(w)= +\infty$, looking at
\[
\lim_{w\rightarrow0^+} w^{- 1+\frac{1}{m} }\int_{]0,+\infty[}
v^{-{1}/m} e^{-( w /v)^{{1}/{\sigma}}}\,d\mu(v).
\]
The last integral, for any $w< 1$, is greater than
\[
\int_{]0,+\infty[} v^{-{1}/m} e^{-(1 /v)^{{1}/{\sigma
}}}\,d\mu(v),
\]
so that it is sufficient to prove that this integral is a positive
constant. If not, this integral would be equal to zero, and this
happens only if the support of $\mu$ is included in ]$-\infty, 0$].
By the martingale connection (\ref{martingaleconnection}), this would
imply that the support of $W^{\mathit{CT}}$ is included in ]$-\infty, 0$],
which is not the case because of Proposition \ref{supportWCT}.

The result on the limit of $p$ at $0^-$ is proved the same way.
Differentiability is immediate by dominated convergence and
monotonicity comes from derivation of
formula (\ref{defp}).
\end{pf}

\begin{Rem}
The distribution of $W^{\mathit{CT}}$ is not symmetric around $0$ (the
expectation equals $\frac b S\not= 0$ when one starts with only one red ball).
\end{Rem}

\subsection{Fourier inversion}\label{inversionFourier}

The characteristic function $\mathcal{F}$ is not integrable.
Nevertheless, formulas (\ref{systemeFourier}) and (\ref{equiv}),
imply straightforwardly that, for any
real $x\not= 0$,\vspace*{-2pt}
\[
\mathcal{F}'(x) = \frac1{mx} \mathcal{F}(x) [ \mathcal
{F}^a(x)\mathcal{G}^b(x) -1
]\vspace*{-2pt}
\]
and that $\mathcal{F}'$ is in ${\rm L}^1$.
Theorem \ref{density2} gives an explicit expression of the density of
$W^{\mathit{CT}}$ by means of inverse
Fourier transform of $\mathcal{F}'$, completing Proposition \ref{density1}.

\begin{Th}
\label{density2}
The density $p$ on $\mathbb{R}$ of the random variable $W^{\mathit{CT}}$ is
given, for
any $x\not= 0$, by\vspace*{-2pt}
%
\begin{equation}
\label{densite}
p(x) = \frac1{2i\pi x} \int_{{\mathbb{R}}} e^{-itx} \mathcal{F}'(t)
\, dt.\vspace*{-2pt}
\end{equation}
\end{Th}

\begin{pf}
Let $F$ be the probability distribution function of $W^{\mathit{CT}}$.
We are going to show that $\forall x\not= 0$,\vspace*{-2pt}
%
\begin{equation}
\label{sufficient}
\lim_{h\rightarrow0} \frac{F(x+h) - F(x)}{h} = \frac1{2i\pi x} \int
_{{\mathbb{R}}} e^{-itx} \mathcal{F}'(t) \, dt,\vspace*{-2pt}
\end{equation}
which is sufficient to prove that $W^{\mathit{CT}}$ admits a continuous density
given by (\ref{densite}). 

For any $h\neq0$, let $d_h$ be the function defined on $\mathbb
{R}\setminus
\{ 0\}$ by
\[
d_h(t):= \frac{1-e^{-ith}}{ith}
\]
and continuated by continuity at $0$.
It follows from the general Fourier inversion theorem
(see, for instance, Lukacs \cite{Lukacs}, Theorem 3.2.1, page~38)
that $\forall x\in\mathbb{R}$, $\forall h\not= 0$, since $x$ and
$x+h$ are
continuity points of
$F$ (remember that $F$ is continuous because its characteristic
function tends to $0$ at infinity),\vspace*{-2pt}
\[
\frac{F(x+h) - F(x)}{h} = \lim_{T\rightarrow+\infty} I_{T,h}(x),\vspace*{-2pt}
\]
where\vspace*{-2pt}
\[
I_{T,h} (x): = \frac1 {2\pi} \int_{-T}^T e^{-itx}d_h(t)\mathcal
{F}(t) \, dt.\vspace*{-2pt}
\]
Integrating by parts implies that, for any $x\neq0$,\vspace*{-2pt}
\[
I_{T,h}(x)= I_{T,h}^{(1)} (x) + I_{T,h}^{(2)}(x) + I_{T,h}^{(3)}(x)\vspace*{-2pt}
\]
where\vspace*{-2pt}
\[
\cases{
I_{T,h}^{(1)} (x)=
\dfrac1{2\pi} \biggl[
-\dfrac{e^{-iTx} }{ix} d_h(T)\mathcal{F}(T)+\dfrac{e^{iTx} }{ix}
d_h(-T)\mathcal{F}(-T)
\biggr] ,
\vspace*{2pt}\cr
I_{T,h}^{(2)}(x)=
\dfrac1{2i\pi x} \int_{-T}^T e^{-itx} d_h(t) \mathcal{F}'(t) \,dt,
\vspace*{2pt}\cr
I_{T,h}^{(3)}(x)=
 \dfrac1{2i\pi x} \int_{-T}^T e^{-itx} d'_h(t) \mathcal{F}(t) \,dt.\vspace*{-1pt}
}
\]
It is elementary to see that $d_h(t) $ has the following properties:
$\forall h \not= 0, \forall t\not= 0$,\vspace*{-2pt}
%
\begin{eqnarray}
\label{P1}
| d_h(t) | &=& \biggl| \frac{\sin{th}/ 2}{ {th}/ 2} \biggr|
\leq\min\biggl\{ 1, \frac2{|th|}\biggr\},
\\[-2pt]
\label{P2}
| d_h'(t) | &\leq&\min\biggl\{ \frac{|h|} 2, \frac2 {|t|}\biggr\}.
\end{eqnarray}
Since $\mathcal{F}$ is bounded (it is a characteristic function) and since
$d_h$ tends to $0$ at infinity,\vspace*{-2pt}
\[
\lim_{T\rightarrow+\infty} I_{T,h}^{(1)} (x) = 0.\vspace*{-2pt}
\]
Since $\mathcal{F}' \in{\rm L}^1$, (\ref{P1}) and Lebesgue dominated
convergence theorem lead to\vspace*{-2pt}
\[
\lim_{T\rightarrow+\infty} I_{T,h}^{(2)} (x) = \frac1{2i\pi x} \int
_{{\mathbb{R}}} e^{-itx} d_h(t) \mathcal{F}'(t)\, dt .\vspace*{-2pt}
\]
At least, (\ref{P2}) implies that $d_h'\mathcal{F} \in{\rm L}^1$ so that,
by dominated convergence,\vspace*{-2pt}
\[
\lim_{T\rightarrow+\infty} I_{T,h}^{(3)} (x)= \frac1{2i\pi x} \int
_{{\mathbb{R}}}e^{-itx} d_h'(t) \mathcal{F}(t)\, dt .\vspace*{-2pt}
\]
So, for any $x\neq0$ and $h\neq0$,\vspace*{-1pt}
\[
\frac{F(x+h) - F(x)}{h} = \frac1{2i\pi x} \int_{{\mathbb{R}}} e^{-itx}d_h(t)
\mathcal{F}'(t)\, dt + \frac1{2i\pi x} \int_{{\mathbb{R}}} e^{-itx}
d_h'(t) \mathcal{F}(t)\, dt .\vspace*{-1pt}
\]
To get (\ref{sufficient}), it is now sufficient to take the limit when
$h\rightarrow0$, using dominated convergence and (\ref{P2}).
\end{pf}

\begin{Rem*}
We have not found the following result in the literature but
the arguments of this proof lead to the following proposition.
\end{Rem*}

\begin{Prop}
Let $\mathcal{F}$ be the characteristic function of a probability
distribution function $F$. 
Suppose that $\mathcal{F}$ is derivable, $\mathcal{F}' \in{\rm L}^1$
($\mathcal{F}$ is not necessarily in~${\rm L}^1$) and $\frac{\mathcal
{F}(t)} t \in
{\rm L}^1 $. Then $F$ admits a density $p$ given for all $x\not=0$ by
\[
p(x) = \frac1{2i\pi x} \int_{{\mathbb{R}}} e^{-itx} \mathcal{F}'(t)
\, dt.
\]
\end{Prop}


\section{Concluding remarks}
\label{conclusion}

\subsection{More colors}

The same questions arise naturally for limit laws of large urn
processes with any finite number of
colors.
Embedding in continuous time, martingale connection, dislocation
equations on elementary limit
distributions and differential system (\ref{systemeFourier}) on
Fourier transforms or on formal Laplace
power series can be generalized.
However, the resolution of (\ref{systemeFourier}) relies on the
question of its integrability, even if
an explicit closed form of its solutions may not be necessary to derive
properties of the corresponding
distributions.

The space requirements of an $m$-ary search tree is a special case of
P\'olya--Eggenberger urn process
with $m-1$ colors (see \cite{ChPo}, for example).
Because of the negativeness of the diagonal entries $-1, -2, \ldots
, -(m-1)$ of its replacement
matrix, the corresponding continuous-time Markov process is not a
branching process.
However, the discrete-time node process of an $m$-ary search tree can
be embedded into a
branching process.
When $m\geq27$, the corresponding limit laws can be studied with the
same method as in the present
paper.
This is the subject of a forthcoming companion paper.

\subsection{Laplace series}
\label{serieLaplace}

Remember from Section \ref{dislocationlimite} that $X$ ({resp.},
$Y$) is the martingale limit
$W^{\mathit{CT}}$ of the continuous-time urn process starting from $(1,0)$
[{resp.}, from $(0,1)$].
For $n\geq0$, let
\[
a_n=\mathbb{E}(X^n)\quad \mbox{and}\quad b_n=\mathbb{E}(Y^n),
\]
and let $F$ and $G$ be the Laplace series of $X$ and $Y$, {\it\textup
{that is,}} the formal exponential series of the moments:
\[
F(T)=\sum_{n\geq0}\frac{a_n}{n!}T^n
\quad \mbox{and}\quad
G(T)=\sum_{n\geq0}\frac{b_n}{n!}T^n\in\mathbb{R}[[T]].
\]
From equations
(\ref{dislocationProjetee}),
we write recursion formulae relating $(a_k)_{0\leq k\leq n}$ and\break
$(b_k)_{0\leq k\leq n}$.
Thanks to the multinomial formula, they arrange themselves into the
differential system with boundary
conditions:
%
\begin{equation}
\label{systemeFormel}
\cases{
F(T)+mTF'(T)=F(T)^{a+1}G(T)^b, \cr
G(T)+mTG'(T)=F(T)^{c}G(T)^{d+1}, \cr
F(0)=G(0)=1, \cr
 F'(0)=\dfrac b S
\quad \mbox{and}\quad
G'(0)=- \dfrac{c} S.
}
\end{equation}
The fact that the urn is large implies that equations (\ref
{systemeFormel}) characterize the
moments of $X$ and $Y$.
Indeed, proceed by recursion:
for any $n\geq2$, $v_n=(a_n,b_n)$ is the solution of a linear system
of the form
$(R-nmI)(v_n)={}$[polynomial function of $v_1,\ldots,v_{n-1}$], $R$
being the replacement matrix of
the process (\ref{matriceUrne}).
Since the urn is large, $nm>nS/2\geq S$ so that $nm$ is not an eigenvalue of $R$.

A remarkable fact, which explains why we have worked with
characteristic functions and not with Laplace transforms, is that, for
nontriangular urns, {\it\textup{that is,}} when $bc\neq0$, series
$F$ and $G$ have a radius of convergence equal to $0$ (Corollary~\ref
{rayonnul}).


\subsection{Question}

The main theorem provides a family of distributions, those of the
$W^{\mathit{CT}}$, indexed by the three
parameters $S,m,b$ of the urn and by the initial condition $(\alpha,
\beta)$.
A challenging question is: can the physical relations between these
distributions be translated into
relations between the Abelian integrals?
In otherwords, can the addition formulas between Abelian integrals be
interpreted by a
combinatorial/probabilistic approach using these distributions?

\section*{Acknowledgments}

The authors warmly thank Philippe Flajolet for stimulating discussions,
Brigitte Chauvin being welcome in Project Algorithms at INRIA Rocquencourt.


%

\printaddresses

\end{document}